\documentclass[12pt]{article}
\usepackage{amsmath}
\usepackage{amssymb,mathrsfs}
\newcommand{\BD}{{\mathbb D}}

\newcommand{\BN}{{\mathbb N}}

\newcommand{\BZ}{{\mathbb Z}}
\newcommand{\sD}{{\cal D}}
\newcommand{\sF}{{\cal F}}
\newcommand{\sH}{{\cal H}}

\newcommand{\sK}{{\cal K}}
\newcommand{\sM}{{\cal M}}

\newcommand{\sU}{{\cal U}}\newcommand{\sV}{{\cal V}}
\newcommand{\sW}{{\cal W}}\newcommand{\sX}{{\cal X}}
\newcommand{\sY}{{\cal Y}}
\newcommand{\eH}{{\mathbf H}}

\newcommand{\eS}{{\mathbf S}}

\newcommand{\Ga}{\Gamma}
\newcommand{\De}{\Delta}

\newcommand{\Th}{\Theta}

\newcommand{\la}{\lambda}\newcommand{\La}{\Lambda}

\newcommand{\om}{\omega}\newcommand{\Om}{\Omega}
\newcommand{\bpr}{{\noindent\textbf{Proof.}\ }}
\newcommand{\epr}{{\hfill $\Box$}\bigskip}
\newcommand{\im}{\textup{Im\,}}
\newcommand{\kr}{\textup{Ker\,}}
\newcommand{\mat}[2]{\ensuremath{\left[\begin{array}{#1}
#2
\end{array} \right]}}
\newcommand{\ov}[1]{{\overline{#1}}}

\newcommand{\inn}[2]{\ensuremath{\langle #1,#2 \rangle}}
\newcommand{\tu}[1]{\textup{#1}}
\newcommand{\half}{\frac{1}{2}}
\newcommand{\ands}{\quad\mbox{and}\quad}
\newcommand{\ons}{\mbox{ on }}
\newcommand{\tilB}{\tilde{B}}

\newcommand{\tilM}{\tilde{M}}

\newcommand{\tilX}{\tilde{X}}

\newcommand{\LDS}{\{A,T',U',R,Q\}}
\newcommand{\Pa}{\Phi}
\newcommand{\Pb}{\tilde\Phi}
\newcommand{\Pc}{\hat\Phi}
\newcommand{\Va}{V}
\newcommand{\Vb}{\tilde V}
\newcommand{\tilGa}{\tilde{\Gamma}}
\newcommand{\Four}{\mathscr{F}}
\newcommand{\Ops}{\mathscr{L}}

\newcommand{\spec}{r_\tu{spec}}
\newcommand{\TA}{T_\tu{state}}
\newcommand{\Lat}{\Lambda^\times}
\newtheorem{theorem}{Theorem}[section]
\newtheorem{corollary}[theorem]{Corollary}

\newtheorem{proposition}[theorem]{Proposition}

\begin{document}
\title{Relaxed commutant lifting and a relaxed Nehari problem:
Redheffer state space formulas}
\author{S. ter Horst}
\date{}
\maketitle
\begin{abstract}
The description of all solutions to the relaxed commutant lifting
problem in terms of an underlying contraction, obtained earlier in
joint work of the author with A.E. Frazho and M.A. Kaashoek, is
transformed into a linear fractional Redheffer state space form.
Under certain additional conditions the coeff{}icient functions in this
representation are described explicitly in terms of the original data.
The main theorem is a generalization of the Redheffer description of
all solutions to the classical commutant lifting problem. To
illustrate the result a relaxed version of the Nehari extension
problem is considered, and an explicit Redheffer description of all
its solutions is given, assuming that a certain truncated Hankel
operator is a strict contraction. The latter result is specif{}ied
further for two special cases.
\end{abstract}

%%%%%%%%%%%%%%%%%%%%%%%%%%%%%%%%%%%%%%%%%%%%%%%%%%%%%%%%%%%%%%%%%%%%%%%%
%%%%%%%%%%%%%%%%%%%%%%%%%%%%%%%%%%%%%%%%%%%%%%%%%%%%%%%%%%%%%%%%%%%%%%%%
%%%%%%%%%%%%%%%%%%%%%%%%%%%%%%%%%%%%%%%%%%%%%%%%%%%%%%%%%%%%%%%%%%%%%%%%
%\setcounter{section}{-1}
\section{Introduction}\label{sec:intro}

The classical commutant lifting theorem, which was obtained by
Sz.-Nagy-Foias \cite{NF68} and originated from the work of Sarason
\cite{S67}, has been used to solve, among other things, a large
number of metric constraint interpolation and extension problems;
see \cite{FFGK98} for a recent overview. In \cite{FFK02a}, extending
the classical theory, a  relaxed commutant lifting problem is
introduced, and a particular (so-called central) solution,
satisfying a maximum entropy condition, is obtained and used to
solve a number of relaxed versions of the classical interpolation
problems. Descriptions of all solutions to the relaxed commutant
lifting problem are given in \cite{FtHK06a}, \cite{LT06} and
\cite{FtHK06b}. The ones in \cite{FtHK06a} and \cite{FtHK06b} are
 in terms of Schur class functions, whereas \cite{LT06} uses a choice
sequences approach. The present paper can be seen as an addition to
\cite{FtHK06b}, where the description of all solutions (see Theorem
\ref{th:FtHK} below) is given in terms of an underlying contraction.
The aim of the present paper is to present a linear fractional
Redheffer type description of all solutions which is explicit in
terms of the original data. For this description some additional
conditions on the data are needed. The main result (see Theorem
\ref{mainth} below) generalizes the corresponding result for the
classical commutant lifting problem given in \cite{FFGK98}.

To be more precise, recall that a {\em lifting data set} is a set
$\Om=\LDS$ consisting of f{}ive Hilbert space operators. The operator
$A$ is a contraction mapping $\sH$ into ${\sH}'$, the operator $U'$
on $\sK'$ is a minimal isometric lifting of the contraction  $T'$ on
$\sH'$, and $R$ and $Q$ are operators from $\sH_0$ to $\sH$,
satisfying the following constraints:
\begin{equation}\label{lds}
T'AR=AQ  \ands R^*R\leq Q^*Q.
\end{equation}
Given a lifting data set $\Om$ as above, the {\em relaxed commutant
lifting problem} is to describe all contractions $B$ from $\sH$ to
$\sK'$ such that
\begin{equation}\label{rclt}
\Pi_{\sH'}B=A\ands U'BR=BQ.
\end{equation}
Here $\Pi_{\sH'}$ is the orthogonal projection from $\sK'$ onto the
subspace $\sH'$. A contraction $B$ from $\sH$ into $\sK'$ satisfying
(\ref{rclt}) will be called a {\em contractive interpolant} for
$\Om$. Hence the relaxed commutant lifting problem is to describe
all contractive interpolants for the lifting data set $\Om$.

Without loss of generality we can, and will, assume that $U'$ is the
Sz.-Nagy-Sch\"affer isometric lifting of $T'$
(\cite{FF90}, Section VI.3), that is, $\sK'$ is the direct sum of
$\sH'$ and the Hardy space $H^2(\sD_{T'})$ and
\begin{equation}\label{szns}
U'=\mat{cc}{T'&0\\E_{\sD_{T'}}D_{T'} & S_{\sD_{T'}}}\ons
\mat{c}{\sH'\\H^2(\sD_{T'})}.
\end{equation}
Here we follow the convention that for a contraction $C$ the symbol
$D_C$ denotes the positive square root of $I-C^*C$, while $\sD_C$
stands for the closure of the range of $D_C$. As usual $D_C$ and
$\sD_C$ are referred to as the {\em defect operator} and {\em defect
space} of $C$, respectively. Moreover, for a Hilbert space $\sY$ the
unilateral shift on the Hardy space $H^2(\sY)$ is denoted by
$S_\sY$, and we write $E_\sY$ for the canonical embedding of $\sY$
onto the subspace of constant functions in $H^2(\sY)$, that is,
$(E_\sY v)(\la)=v$ for all $\la\in\BD$ and each $v\in\sY$.

Put
\begin{equation}\label{Dcirc}
D_\circ=(Q^*Q-R^*R)^\half\ands\sD_\circ=\ov{D_\circ\sH_0}.
\end{equation}
Then (\ref{lds}) implies that
\begin{equation}\label{equal}
Q^*D_A^2Q=D_\circ^2+R^*A^*D_{T'}^2AR+R^*D_A^2R.
\end{equation}
Indeed, this is the case because for all $h\in\sH_0$ we have
\begin{eqnarray*}
\|D_AQh\|^2&=&\|Qh\|^2-\|AQh\|^2
=\|D_\circ h\|^2+\|Rh\|^2-\|T'ARh\|^2\\
&=&\|D_\circ h\|^2+\|ARh\|^2-\|T'ARh\|^2+\|Rh\|^2-\|ARh\|^2\\
&=&\|D_\circ h\|^2+\|D_{T'}ARh\|^2+\|D_ARh\|^2.
\end{eqnarray*}
With the lifting data set $\Om$ we associate a contraction $\om$
def{}ined by
\begin{equation}\label{om}
\om:\sF\to\sD_{T'}\oplus\sD_A,\quad\om
D_AQ=\mat{cc}{D_{T'}AR\\D_AR}, \quad\sF=\ov{D_AQ\sH_0}.
\end{equation}
The relation in (\ref{equal}) guarantees that $\om$ is contractive.
Moreover, $\om$ is an isometry if and only if $D_\circ=0$. In terms of
the operators from the lifting data set this is equivalent to
$R^*R=Q^*Q$. The f{}irst main theorem from \cite{FtHK06b} can now be
formulated as follows.

%%%%%%%%%%%%%%%%%%%%%%%%%%%%%%%%%%%
\begin{theorem}\label{th:FtHK}
Let $\Om=\LDS$ be a lifting data set with $U'$ the
Sz.-Nagy-Sch\"affer isometric lifting of $T'$, and let $B$ be an
operator from $\sH$ into $\sH'\oplus H^2(\sD_{T'})$. Then $B$ is a
contractive interpolant for $\Om$ if and only if $B$ admits a
representation of the form
\begin{equation}\label{solsA}
B=\mat{c}{A\\\tilGa D_A}:\sH\to\mat{c}{\sH'\\H^2(\sD_{T'})},
\end{equation}
where $\tilGa$ is the operator mapping $\sD_A$ into $H^2(\sD_{T'})$
given by
\begin{equation}\label{solsB}
(\tilGa d)(\la)=\Pi_{\sD_{T'}}Z(\la)(I-\la\Pi_{\sD_A}Z(\la))^{-1}d
\quad\quad(d\in\sD_A,\la\in\BD)
\end{equation}
for some Schur class function $Z$ from
$\eS(\sD_A,\sD_{T'}\oplus\sD_A)$ with $Z(\la)|\sF=\om$ for each
$\la\in\BD$.
\end{theorem}

Let us explain in some detail the notations that are used in the above
theorem or will appear in the sequel. Throughout capital calligraphic
letters denote Hilbert spaces. The Hilbert space direct sum of $\sU$
and $\sY$ is denoted by
\[
\sU\oplus \sY \quad \mbox{or by}\quad \left[
\begin{array}{c}
\sU\\ \sY
\end{array}
\right].
\]
An {\em operator} is a bounded linear transformation acting between
Hilbert spaces. With $\Ops(\sU, \sY)$ we denote the set of all
operators from $\sU$ into $\sY$. The {\em identity operator} on the
space $\sU$ is denoted by $I_{\sU}$, or just by $I$ when the
underlying space is clear from the context. By def{}inition, a {\em
subspace} is a closed linear manifold. Let $\sM$ be a subspace of
$\sU$. Then $\sU\ominus \sM$ stands for the {\em orthogonal
complement} of $\sM$ in $\sU$. We follow the convention that the symbol
$\Pi_\sM$ denotes the {\em orthogonal projection} from $\sU$ onto $\sM$
viewed as an operator from $\sU$ to $\sM$, whereas $P_\sM$ stands for
the {\em orthogonal projection} from $\sU$ onto $\sM$ acting as an
operator on $\sU$. Note that with this notation $\Pi_\sM^*$ is the
{\em canonical embedding} of $\sM$ into $\sU$ and
$P_\sM=\Pi_\sM^*\Pi_\sM$. An operator $C$ on $\sU$
is referred to as positive def{}inite (notation: $C>0$) if $C$ is
invertible and positive (i.e., $\inn{Cu}{u}\geq0$ for each $u\in\sU$).
An operator $N$ in $\Ops(\sU,\sY)$ is called {\em left invertible} if
there exists an operator $M$ in $\Ops(\sY,\sU)$ such that $MN=I_\sU$.
In that case $M$ is called a {\em left inverse} of $N$. Note that $N$
is left invertible if and only if $N^*N$ is positive def{}inite.
By def{}inition, the symbol $\eH^\infty(\sU,\sY)$ stands for the class
of uniformly bounded analytic functions on the open unit disc $\BD$
with values in $\Ops(\sU,\sY)$. The set $\eH^\infty(\sU,\sY)$ forms a
Banach space with respect to the supremum norm. In that case
$\eS(\sU,\sY)$ denotes the closed unit disc in $\eH^\infty(\sU,\sY)$, and is referred to as the
{\em Schur class} associated with $\sU$ and $\sY$. Functions from
$\eS(\sU,\sY)$ are called {\em Schur class functions}. Finally, by
$\eH^2(\sU,\sY)$ we denote the class of all functions $F$ on $\BD$ with
values in $\Ops(\sU,\sY)$ such that for each $u$ in $\sU$ the map
$\la\mapsto F(\la)u$ def{}ines a function in the Hardy class $H^2(\sY)$.
Such a function is automatically analytic on $\BD$. Obviously
$\eH^\infty(\sU,\sY)$ is properly contained in $\eH^2(\sU,\sY)$.

The next theorem is the main result of the present paper.

%%%%%%%%%%%%%%%%%%%%%%%%%%%%%%%%%%%
\begin{theorem}\label{mainth}
Let $\Om=\LDS$ be a lifting data set with $U'$ the
Sz.-Nagy-Sch\"affer isometric lifting of $T'$. Assume that $A$ is
a strict contraction and $R$ has a left inverse. Then an operator
$B$ mapping $\sH$ into $\sH'\oplus H^2(\sD_{T'})$ is a contractive
interpolant for $\Om$ if and only if $B$ admits a representation of
the form
\begin{equation}\label{sols1}
B=\mat{c}{A\\\Ga}:\sH\to\mat{c}{\sH'\\H^2(\sD_{T'})},
\end{equation}
where $\Ga$ is the operator mapping $\sH$ into $H^2(\sD_{T'})$ given by
\begin{equation}\label{sols2}
(\Ga h)(\la)=\Pa_{22}(\la)h+\Pa_{21}(\la)\Va(\la)
(I-\Pa_{11}(\la)\Va(\la))^{-1}\Pa_{12}(\la)h\quad(h\in\sH,\la\in\BD)
\end{equation}
for some Schur class function $\Va$ from $\eS(\kr
Q^*,\sD_\circ\oplus\sD_{T'}\oplus\kr R^*)$. Moreover, any such function $V$
def{}ines an operator $\Ga$ mapping $\sH$ into $H^2(\sD_{T'})$ via
\tu{(\ref{sols2})}. Here $\Pa_{11}$ and $\Pa_{21}$ are
functions from the Schur classes $\eS(\sD_\circ\oplus\sD_{T'}\oplus\kr R^*,\kr Q^*)$
and $\eS(\sD_\circ\oplus\sD_{T'}\oplus\kr R^*,\sD_{T'})$,
respectively, and $\Pa_{12}$ and $\Pa_{22}$ are functions from
$\eH^2(\sH,\kr Q^*)$ and $\eH^2(\sH,\sD_{T'})$, respectively, and
these functions are given by
\begin{equation}\label{sols3}
\begin{array}{rcl}
\Pa_{11}(\la)&=&\la X_3(I-\la X_1)^{-1}X_2,\\[.1cm]
\Pa_{12}(\la)&=&X_3(I-\la X_1)^{-1},\\[.1cm]
\Pa_{21}(\la)&=&X_5+\la X_4(I-\la X_1)^{-1}X_2,\\[.1cm]
\Pa_{22}(\la)&=&X_4(I-\la X_1)^{-1},
\end{array}
\quad\quad\quad\quad\quad\quad\quad(\la\in\BD)
\end{equation}
where $X_1,\ldots,X_5$ are the operators def{}ined by
\begin{equation}\label{X1X5}
\begin{array}{rcl}
X_1&=&R(Q^*D_A^2Q)^{-1}Q^*D_A^2\ons\sH,\\[.1cm]
X_2&=&-R(R^*D_A^2R)^{-1}J^*\De_\Om^{-\half}
\Pi_{\sD_\circ\oplus\sD_{T'}}+\\[.1cm]
&&\quad\quad\quad-D_A^{-2}\Pi_{\kr R^*}^*\De_R^{-\half}\Pi_{\kr R^*}
:\sD_\circ\oplus\sD_{T'}\oplus\kr R^*\to\sH,\\[.1cm]
X_3&=&\De_Q^{-\half}\Pi_{\kr Q^*}:\sH\to\kr Q^*\\[.1cm]
X_4&=&D_{T'}AR(Q^*D_A^2Q)^{-1}Q^*D_A^2:\sH\to\sD_{T'},\\[.1cm]
X_5&=&\Pi_{\sD_{T'}}\De_\Om^{-\half}\Pi_{\sD_\circ\oplus\sD_{T'}}
:\sD_\circ\oplus\sD_{T'}\oplus\kr R^*\to\sD_{T'},
\end{array}
\end{equation}
and $\De_Q$ on $\kr Q^*$, $\De_R$ on $\kr R^*$ and $\De_\Om$ on
$\sD_\circ\oplus\sD_{T'}$ are the positive def{}inite operators
def{}ined by
\begin{equation}\label{Deltas}
\begin{array}{l}
\De_Q=\Pi_{\kr Q^*}D_A^{-2}\Pi_{\kr Q^*}^*,\quad\quad
\De_R=\Pi_{\kr R^*}D_A^{-2}\Pi_{\kr R^*}^*,\\[.1cm]
\De_\Om=I+J(R^*D_A^2R)^{-1}J^*,\quad\mbox{where}\quad
J=\mat{c}{D_\circ\\D_{T'}AR}:\sH_0\to\mat{c}{\sD_\circ\\\sD_{T'}}.
\end{array}
\end{equation}
\end{theorem}

The condition in the above theorem that $A$ is a strict contraction
is equivalent to the requirement that the defect operator $D_A$ of $A$
is positive def{}inite on $\sH$. Moreover, if $R$ is left invertible,
then the second condition in (\ref{lds}) implies that $Q$ is also left
invertible. Thus the combination of both conditions  results in $D_AQ$
and $D_AR$ being left invertible, or equivalently, $Q^*D_A^2Q$ and
$R^*D_A^2R$ being positive def{}inite, both on $\sH_0$. So the extra
assumptions on the lifting data set in Theorem \ref{mainth} imply that
the operators $D_A^{-2}$, $(R^*D_A^2R)^{-1}$ and $(Q^*D_A^2Q)^{-1}$
appearing in Theorem \ref{mainth} are well def{}ined.

When taking the zero function for $\Va$ in Theorem \ref{mainth} we see
that (\ref{sols2}) reduces to $(\Ga h)(\la)=\Pa_{22}(\la)h$ for
$h\in\sH$ and $\la\in\BD$. The solution obtained in this way is
precisely the central solution given in \cite{FFK02a}.

In \cite{FFK02a}, Proposition 5.3, it was shown that the spectral
radius of $X_1$ in (\ref{X1X5}) is strictly less then one if $R$ and
$Q$ are such that $R-\la Q$ is left invertible for each $\la\in\BD$.
Note that in this case the functions $\Pa_{12}$ and $\Pa_{22}$ are
uniformly bounded on $\BD$, that is, $\Pa_{12}$ and $\Pa_{22}$ are
functions from the classes $\eH^\infty(\sH,\kr Q^*)$ and
$\eH^\infty(\sH,\sD_{T'})$, respectively. This remark will be useful
later on when we consider the relaxed Nehari extension problem.

The proof of Theorem \ref{mainth} is given in Section \ref{sec:RED}.
In this section we also derive additional properties of the operator
valued functions in (\ref{sols3}) (see Corollary \ref{cor:multop}
below).

Recall that the classical commutant lifting problem appears when in
the lifting data set $\LDS$ the operator $Q$ is an isometry and $R$
the identity operator on $\sH$, and thus, in particular, $\sH_0=\sH$.
Under these additional conditions Theorem \ref{mainth} reduces to the
f{}irst part of Theorem VI.6.1 in \cite{FFGK98} (see Corollary
\ref{cor:CCL} below for further details). Moreover, in that case
(\ref{sols2}) provides a proper parameterization, that is, there
exists a unique Schur class function $V$ such that $B$ is given by
(\ref{sols1}) and (\ref{sols2}). In general, for an arbitrary lifting
data set formula (\ref{sols2}) does not provide a proper
parameterization. This follows from Theorem 1.2 in \cite{FtHK06b}.

To illustrate Theorem \ref{mainth} we consider a relaxed version of
the operator-valued Nehari extension problem. Let $N$ be a positive
integer $(N>0)$ and let $F_{-1},F_{-2},\ldots$ be a sequence of
operators from $\sU$ to $\sY$ satisfying
$\sum_{n=1}^\infty\|F_{-n}u\|^2<\infty$ for each $u\in\sU$. The {\em
relaxed Nehari extension problem} considered in this paper is to
f{}ind all sequences of operators $H_0,H_1,\ldots$ from $\sU$ to $\sY$
with the property that $\sum_{n=0}^\infty\|H_nu\|^2<\infty$ for each
$u\in\sU$, and such that the operator from $\sU^N$ into $\ell^2(\sY)$
given by the operator matrix representation
\begin{equation}\label{L}
\mat{cccc}
{\vdots&\vdots&&\vdots\\
F_{-2}&F_{-3}&\cdots&F_{-(N+1)}\\
F_{-1}&F_{-2}&\cdots&F_{-N}\\
\fbox{$H_0$}
&F_{-1}&\cdots&F_{-(N-1)}\\
H_1&H_0&\ddots&\vdots\\
\vdots&\ddots&\ddots&F_{-1}\\
\vdots&&\ddots&H_0\\
\vdots&&&\vdots\\}:\sU^N\to\ell^2(\sY)
\end{equation}
has operator norm at most one. Here $\sU^N$ is the Hilbert space
direct sum of $N$ copies of $\sU$, and $\ell^2(\sY)$ is the Hilbert
space of bilateral square summable sequences $(y_n)_{n\in\BZ}$ with
entries in $\sY$. As usual $\BZ$ stands for the set of all integers.
The box in (\ref{L}) indicates the zero position in $\ell^2(\sY)$. A
sequence of operators $(H_n)_{n\in\BN}$ from $\sU$ to $\sY$ that
forms a solution to the relaxed Nehari problem is referred to as an
{\em $N$-complementary sequence} associated with $(F_{-n-1})_{n\in\BN}$,
or just an $N$-complementary sequence if no confusion concerning the
sequence $(F_{-n-1})_{n\in\BN}$ can arise. Here $\BN$ stands for the
set of nonnegative integers (i.e., with zero included). The setup for
this problem resembles the way relaxed versions of the Schur,
Nevanlinna-Pick and Sarason interpolation problems where formulated in
\cite{FFK02a}. For this relaxed Nehari problem to be solvable it is
necessary that the operator given by
\begin{equation}\label{NHankel}
\mat{cccc} {\vdots&\vdots&&\vdots\\
F_{-3}&F_{-4}&\cdots&F_{-(N+2)}\\
F_{-2}&F_{-3}&\cdots&F_{-(N+1)}\\
F_{-1}&F_{-2}&\cdots&F_{-N}}:\sU^N\to\ell^2_-(\sY)
\end{equation}
is a contraction. Here  $\ell^2_-(\sY)$ stands for the Hilbert space of all
square summable sequences $(\ldots,y_{-2},y_{-1})$ with entries in
$\sY$.

In Section \ref{sec:Nehari} the relaxed Nehari extension problem is
put into a relaxed commutant lifting setting, which yields that the
condition that the operator in (\ref{NHankel}) is a contraction is
not only necessary but also suff{}icient. The latter can also be seen
by repeatedly applying Parrott's lemma (see Corollary IV.3.6 in
\cite{FF90}). We use Theorem \ref{mainth} to give a Redheffer
description (in Theorem \ref{th:Neh} below) of all $N$-complementary
sequences, under the additional assumption that the operator in
(\ref{NHankel}) is a strict contraction. In addition, we specify
Theorem \ref{th:Neh} for two special cases, namely, when $N=1$
(see Corollary \ref{cor:N=1} below), and when $F_n=0$ for
$n=-1,-2,\ldots$ (see Corollary \ref{cor:Fn=0} below).

%%%%%%%%%%%%%%%%%%%%%%%%%%%%%%%%%%%%%%%%%%%%%%%%%%%%%%%%%%%%%%%%%%%%%%%%
%%%%%%%%%%%%%%%%%%%%%%%%%%%%%%%%%%%%%%%%%%%%%%%%%%%%%%%%%%%%%%%%%%%%%%%%
\section{Redheffer representations and proof of the main theorem}
\label{sec:RED}%\setcounter{equation}{0}

The aim in this section is to prove Theorem \ref{mainth}. The main
tool to achieve this objective is Proposition \ref{pr:RED} given in
the next subsection. This proposition gives a general scheme for
rewriting the description in (\ref{solsB}) into one of the type
(\ref{sols2}).

%%%%%%%%%%%%%%%%%%%%%%%%%%%%%%%%%%%%%%%%%%%%%%%%%%%%%%%%%%%%%%%%%%%%%%%%
\subsection{Redheffer type descriptions}

First we introduce some notation. Let $C$ be an operator on $\sU$.
Then $C$ is said to be pointwise stable if $C^nu$ converges to zero
as $n$ goes to inf{}inity, for each $u\in\sU$. Recall that $C$ is
pointwise stable if the spectral radius of $C$, denoted by
$\spec(C)$, is strictly less then~1. Next, f{}ix a function $H$  from
$\eH^2(\sU,\sY)$. With $H$ we associate an operator $\Ga_H$ mapping
$\sU$ into $H^2(\sY)$ def{}ined by
\begin{equation}\label{multH2}
(\Ga_H u)(\la)=H(\la)u\quad\quad(\la\in\BD,u\in\sU).
\end{equation}
On the other hand, if $\Ga$ is any operator from $\sU$ into
$H^2(\sY)$, then there exists a unique function $G$ from
$\eH^2(\sU,\sY)$ such that $\Ga=\Ga_G$. This function is given by
$G(\la)u=(\Ga u)(\la)$ for $u\in\sU$ and $\la\in\BD$. Now assume, in
addition, that $H$ is from $\eH^\infty(\sU,\sY)$. We associate with
$H$, in the usual way, a multiplication operator $M_H$  mapping
$H^2(\sU)$ into $H^2(\sY)$ def{}ined by
\begin{equation*}\label{multS}
(M_Hf)(\la)=H(\la)f(\la)\quad\quad(\la\in\BD,f\in H^2(\sU)).
\end{equation*}
The operator $M_H$ is called the {\em multiplication operator}
def{}ined by $H$ and its norm is given by
\begin{equation*}\label{multnorm}
\|M_H\|=\sup_{\la\in\BD}\|H(\la)\|.
\end{equation*}
In particular, $M_H$ is a contraction if and only if $H$ is from the
Schur class $\eS(\sU,\sY)$.

%%%%%%%%%%%%%%%%%%%%%%%%%%%%%%%%%%%
\begin{proposition}\label{pr:RED}
Assume that $Z$ is a Schur class function from
$\eS(\sU,\sY\oplus\sU)$ given by
\begin{equation}\label{Somfuncs}
Z(\la)
=\mat{c}{\tilX_4\\\tilX_1}+\mat{c}{\tilX_5\\\tilX_2}\Va(\la)\tilX_3
:\sU\to\mat{c}{\sY\\\sU}\quad\quad(\la\in\BD),
\end{equation}
where $V$ is from the Schur class $\eS(\sV,\sW)$ and
\begin{equation}\label{tilX}
\tilX=\mat{cc}{\tilX_1&\tilX_2\\\tilX_3&0\\\tilX_4&\tilX_5}
:\mat{c}{\sU\\\sW}\to\mat{c}{\sU\\\sV\\\sY}\mbox{ is a contraction.}
\end{equation}
Def{}ine operator-valued functions
$\Pb_{11}$, $\Pb_{12}$, $\Pb_{21}$ and $\Pb_{22}$ by
\begin{equation}\label{Pb}
\begin{array}{rcl}
\Pb_{11}(\la)&=&\la\tilX_3(I-\la\tilX_1)^{-1}\tilX_2,\\[.2cm]
\Pb_{12}(\la)&=&\tilX_3(I-\la\tilX_1)^{-1},\\[.2cm]
\Pb_{21}(\la)&=&\tilX_5+\la\tilX_4(I-\la\tilX_1)^{-1}\tilX_2,\\[.2cm]
\Pb_{22}(\la)&=&\tilX_4(I-\la\tilX_1)^{-1}.
\end{array}
\quad\quad\quad\quad\quad\quad\quad(\la\in\BD)
\end{equation}
Then $\Pb_{11}\in\eS(\sW,\sV)$, $\Pb_{21}\in\eS(\sW,\sD_{T'})$,
$\Pb_{12}\in\eH^2(\sU,\sV)$, $\Pb_{22}\in\eH^2(\sU,\sY)$ and
\begin{equation}\label{transform}
\Pi_\sY Z(\la)(I-\la\Pi_\sU Z(\la))^{-1}
=\Pb_{22}(\la)+\Pb_{21}(\la)\Va(\la)(I-\Pb_{11}(\la)\Va(\la))^{-1}
\Pb_{12}(\la)\quad(\la\in\BD).
\end{equation}
Moreover, the operator
\begin{equation}\label{tilmultop}
\tilM=\mat{cc}{M_{\Pb_{11}}&\Ga_{\Pb_{12}}\\M_{\Pb_{21}}&\Ga_{\Pb_{22}}}
:\mat{c}{H^2(\sW)\\\sU}\to \mat{c}{H^2(\sV)\\H^2(\sY)}
\end{equation}
is a contraction which is unitary whenever $\tilX$ in
\tu{(\ref{tilX})} is unitary and $\tilX_1$ is pointwise stable.
\end{proposition}

The formula on the right hand side of (\ref{transform}) is referred to
as a linear fractional Redheffer description. The term Redheffer comes
from scattering theory, see Chapter XIV in \cite{FF90}. Indeed, let
$\Pb_{11}$, $\Pb_{12}$, $\Pb_{21}$ and $\Pb_{22}$ be the functions in
(\ref{Pb}) where $\tilX$ in (\ref{tilX}) is a contraction. Consider
the Redheffer scattering system
\begin{equation}\label{REDsystem}
\mat{c}{g\\y}=\mat{cc}{M_{\Pb_{11}}&\Ga_{\Pb_{12}}\\
M_{\Pb_{21}}&\Ga_{\Pb_{22}}}
\mat{c}{x\\u},\quad
\mbox{subject to }x=M_Vg,
\end{equation}
where $V$ is a Schur class function from $\eS(\sV,\sW)$. Here $u$ is an
element from $\sU$ and the vectors $y$, $g$ and $x$ are functions from
the Hardy spaces $H^2(\sY)$, $H^2(\sV)$ and $H^2(\sW)$, respectively.
Solving (\ref{REDsystem}) we obtain that $y=\tilGa u$, where $\tilGa$
is the operator from $\sU$ into $H^2(\sY)$ def{}ined by the function in
the right hand side of (\ref{transform}), that is,
\[
(\tilGa u)(\la)=\Pb_{22}(\la)u
+\Pb_{21}(\la)\Va(\la)(I-\Pb_{11}(\la)\Va(\la))^{-1}\Pb_{12}(\la)u
\quad\quad(u\in\sU,\la\in\BD).
\]

\noindent{\bf Proof of Proposition \ref{pr:RED}}
Since $\tilX$ is contractive, and thus $\tilX_1$ is
contractive, the functions in (\ref{Pb}) are properly def{}ined, and
analytic on $\BD$. Moreover, we have
\[
\Pi_\sY Z(\la)=\tilX_4+\tilX_5\Va(\la)\tilX_3\ands
\Pi_\sU Z(\la)=\tilX_1+\tilX_2\Va(\la)\tilX_3
\quad\quad(\la\in\BD).
\]
For each $\la\in\BD$ we then obtain that
\begin{eqnarray*}
\tilX_3(I-\la \Pi_\sU Z(\la))^{-1}
&=&\tilX_3(I-\la\tilX_1-\la\tilX_2\Va(\la)\tilX_3)^{-1}\\
&=&\tilX_3(I-\la(I-\la\tilX_1)^{-1}\tilX_2\Va(\la)\tilX_3)^{-1}
(I-\la\tilX_1)^{-1}\\
&=&(I-\la\tilX_3(I-\la\tilX_1)^{-1}\tilX_2\Va(\la))^{-1}
\tilX_3(I-\la\tilX_1)^{-1}\\
&=&(I-\Pb_{11}(\la)\Va(\la))^{-1}\Pb_{12}(\la),
\end{eqnarray*}
and
\begin{eqnarray*}
\tilX_4(I-\la \Pi_\sU Z(\la))^{-1}
&=&\tilX_4(I-\la(I-\la\tilX_1)^{-1}\tilX_2\Va(\la)\tilX_3)^{-1}
(I-\la\tilX_1)^{-1}\\
&=&\tilX_4(I-\la\tilX_1)^{-1}+
\la\tilX_4(I-\la\tilX_1)^{-1}\tilX_2\Va(\la)\tilX_3
(I-\la \Pi_\sU Z(\la))^{-1}\\
&=&\Pb_{22}(\la)+(\Pb_{21}(\la)-\tilX_5)V(\la)
(I-\Pb_{11}(\la)\Va(\la))^{-1}\Pb_{12}(\la).
\end{eqnarray*}
The combination of these two results gives
\begin{eqnarray*}
\Pi_\sY Z(\la)(I-\la \Pi_\sU Z(\la))^{-1}
&=&(\tilX_4+\tilX_5\Va(\la)\tilX_3)(I-\la \Pi_\sU Z(\la))^{-1}\\
&=&\Pb_{22}(\la)+(\Pb_{21}(\la)-\tilX_5+\tilX_5)V(\la)
(I-\Pb_{11}(\la)\Va(\la))^{-1}\Pb_{12}(\la)\\
&=&\Pb_{22}(\la)+\Pb_{21}(\la)V(\la)
(I-\Pb_{11}(\la)\Va(\la))^{-1}\Pb_{12}(\la).
\end{eqnarray*}
So (\ref{transform}) holds.

For the remainder of the proof we use some results from system
theory. The terminology corresponds to that in \cite{FFGK98}. A
contractive system is a quadruple $\Th=\{Z,B,C,D\}$, consisting of
operators $Z$ on a Hilbert space $\sX$, $B$ from $\sU$ to $\sX$, $C$
from $\sX$ to $\sY$ and $D$ mapping $\sU$ into $\sY$
such that the operator matrix
\begin{equation}\label{sysmat}
K_\Th=\mat{cc}{Z&B\\C&D}:\mat{c}{\sX\\\sU}\to\mat{c}{\sX\\\sY}
\mbox{ is a contraction.}
\end{equation}
Let $\Th=\{Z,B,C,D\}$ be a contractive system. Since $Z$ is
contractive, we can def{}ine operator-valued functions $F_\Th$ and
$G_\Th$ on $\BD$ by
\begin{equation}\label{multobs}
F_\Th(\la)=D+\la C(I-\la Z)^{-1}B\ands
G_\Th(\la)=C(I-\la Z)^{-1}\quad\quad(\la\in\BD).
\end{equation}
Here $F_\Th$ is referred to as the transfer function for $\Th$. {}From
the fact that $K_\Th$ in (\ref{sysmat}) is contractive it follows that
$F_\Th\in\eS(\sU,\sY)$ and $G_\Th\in\eH^2(\sX,\sY)$. The operator
$\Ga_{G_\Th}$ from $\sU$ to $H^2(\sY)$ is referred to as the
observability operator for $\Th$. Moreover, we have that
\begin{equation}\label{muop}
\mat{cc}{M_{F_\Th}&\Ga_{G_\Th}}:\mat{c}{H^2(\sU)\\\sX}\to H^2(\sY)
\end{equation}
is a contractive operator which is unitary whenever $K_\Th$ is unitary and
$Z$ is pointwise stable. The statement for the case that $K_\Th$ is
unitary and $Z$ pointwise stable is obtained from Theorem III.10.4 in
\cite{FFGK98}, the statement that the operator in (\ref{muop}) is
contractive in the general case can easily be derived from Theorem
III.10.1 in \cite{FFGK98}, it also follows by specifying the result
from Proposition 1.7.2 in \cite{P99} concerning time-variant systems
for the time-invariant case.

Now put
\begin{equation}\label{sysproof}
\Th=\left\{\tilX_1,\tilX_2,\mat{c}{\tilX_3\\\tilX_4},
\mat{c}{0\\\tilX_5}\right\}.
\end{equation}
Then the operator $K_\Th$ in (\ref{sysmat}) is equal to $\tilX$ in
(\ref{tilX}). So $\Th$ is a contractive system. Moreover, the functions
$F_\Th$ and $G_\Th$ in (\ref{multobs}) are given by
\[
F_\Th=\mat{c}{\Pb_{11}\\\Pb_{21}}\ands G_\Th=\mat{c}{\Pb_{12}\\\
\Pb_{22}}.
\]
The proposition then follows from the theory concerning contractive
systems summed up above and the observation that
\[
\tilM=\mat{cc}{M_{F_\Th}&\Ga_{G_\Th}},\mbox{ where $\Th$ is given by
(\ref{sysproof})}.
\]
Here we identify $H^2(\sV\oplus\sY)$ with $H^2(\sV)\oplus H^2(\sY)$.
\epr

%%%%%%%%%%%%%%%%%%%%%%%%%%%%%%%%%%%%%%%%%%%%%%%%%%%%%%%%%%%%%%%%%%%%%%%%
\subsection{Proof of Theorem \ref{mainth}}

Proposition \ref{pr:RED} suggests that in order to prove Theorem
\ref{mainth} it suffices to show that a Schur class function $Z$ from
$\eS(\sD_A,\sD_{T'}\oplus\sD_A)$ with the property that
$Z(\la)|\sF=\om$ for each $\la\in\BD$ can be expressed as in
(\ref{Somfuncs}), with $\tilX_1$, $\tilX_2$, $\tilX_3$, $\tilX_4$ and $\tilX_5$ the appropriate operators. This is done in the next
proposition under the additional assumption that $A$ is a strict
contraction and $R$ is left invertible.

%%%%%%%%%%%%%%%%%%%%%%%%%%%%%%%%%%%
\begin{proposition}\label{pr:BTR}
Let $\Om=\LDS$ be a lifting data set, and let $Z$ be a Schur class function
from $\eS(\sD_A,\sD_{T'}\oplus\sD_A)$. Assume that $A$ is a strict
contraction and $R$ is left invertible. Then $Z(\la)|\sF=\om$ for
each $\la\in\BD$ if and only if there exists a function $\Va$ from
the Schur class $\eS(\kr Q^*,\sD_\circ\oplus\sD_{T'}\oplus\kr R^*)$
such that
\begin{equation}\label{ZV}
Z(\la)=\mat{c}{X_4D_A^{-1}\\D_AX_1D_A^{-1}}+\mat{c}{X_5\\D_AX_2}
\Va(\la)X_3D_A^{-1}\quad\quad(\la\in\BD),
\end{equation}
where $X_1$, $X_2$, $X_3$, $X_4$ and $X_5$ are the operators def{}ined in
\tu{(\ref{X1X5})}. Finally, $Z$ and $\Va$ def{}ine each other uniquely in
\tu{(\ref{ZV})}.
\end{proposition}

\epr
In \cite{FtHK06b} it was already shown that a function $Z$ from
$\eS(\sD_A,\sD_{T'}\oplus\sD_A)$ satisf{}ies the equality
$Z(\la)|\sF=\om$ for each $\la\in\BD$ if and only if there exists a
Schur class function $\Vb$ from $\eS(\kr Q^*D_A,\sD_{\om^*})$ such that
\begin{equation}\label{ZVb}
Z(\la)=\om\Pi_\sF+D_{\om^*}\Vb(\la)\Pi_{\kr Q^*D_A}
\quad\quad(\la\in\BD),
\end{equation}
where $\om$ is the contraction def{}ined in (\ref{om}).
Moreover, it was shown there that $Z$ and $\Vb$ in (\ref{ZVb})
def{}ine each other uniquely. Recall that $\|A\|<1$ and $R$ being left
invertible imply that $D_AQ$ and $D_AR$ are left invertible, or
equivalently, that $Q^*D_A^2Q$ and $R^*D_A^2R$ are positive
def{}inite on $\sH_0$. Then left inverses of $D_AQ$ and $D_AR$ are given
by
\begin{equation}\label{DAQRinvs}
L_{D_AQ}=(Q^*D_AQ)^{-1}Q^*D_A\ands
L_{D_AR}=(R^*D_AR)^{-1}R^*D_A,
\end{equation}
respectively. Note that $L_{D_AQ}|\kr Q^*D_A=0$ and
$L_{D_AR}|\kr R^*D_A=0$.

The proof of Proposition \ref{pr:BTR} consists of four parts. In Part 1
we show that
\begin{equation*}\label{omform}
\om\Pi_\sF=\mat{c}{X_4\\D_AX_1}D_A^{-1}.
\end{equation*}
Put
\begin{equation*}\label{Y}
Y=\mat{cc}{\De_\Om^{-\half}\Pi_{\sD_{T'}}^*&-\De_\Om^{-\half}JL_{D_AR}\\
0&-\Pi_{\kr R^*D_A}}
:\mat{c}{\sD_{T'}\\\sD_A}\to\mat{c}{\sD_\circ\oplus\sD_{T'}\\\kr
R^*D_A},
\end{equation*}
where $\De_\Om$ and $J$ are the operators in (\ref{Deltas}).
The second part is used to prove that $Y^*Y=D_{\om^*}^2$ and
$\kr Y^*=\{0\}$. In Part 3 we show that
$\im D_A^{-1}\Pi_{\kr Q^*}=\kr Q^* D_A$ and $\im D_A^{-1}\Pi_{\kr
R^*}^*=\kr R^* D_A$.
With these identities we obtain that
\begin{equation}\label{projects}
P_{\kr Q^* D_A}=D_A^{-1}\Pi_{\kr Q^*}^*\De_Q^{-1}\Pi_{\kr
Q^*}D_A^{-1}\mbox{ and } P_{\kr R^* D_A}=D_A^{-1}\Pi_{\kr
R^*}^*\De_R^{-1}\Pi_{\kr R^*}D_A^{-1},
\end{equation}
where $\De_Q$ and $\De_R$ are the operators in (\ref{Deltas}).
The results from the Parts 1 to 3 are then combined in Part 4 to
complete the proof of Proposition \ref{pr:BTR}.\bigskip

\noindent{\bf Part 1.}
Since $\sD_A\ominus\sF=\kr Q^*D_A$, we have
$L_{D_AQ}|\sD_A\ominus\sF=0$, where $L_{D_AQ}$ is the left inverse of
$D_AQ$ given by (\ref{DAQRinvs}). {}From the def{}inition of $\om$ in
(\ref{om}) we then obtain that
\begin{equation}\label{omtrans}
\om\Pi_\sF=\mat{c}{D_{T'}AR\\D_AR}L_{D_AQ}
=\mat{c}{D_{T'}AR(Q^*D_A^2Q)^{-1}Q^*D_A\\D_AR(Q^*D_A^2Q)^{-1}Q^*D_A}
=\mat{c}{X_4\\D_AX_1}D_A^{-1}.
\end{equation}

\noindent{\bf Part 2.}
We f{}irst show that $J^*$ intertwines the
operator $\De_\Om$ with $Q^*D_A^2Q(R^*D_A^2R)^{-1}$, that is,
\begin{equation}\label{K*DeOm}
J^*\De_\Om=Q^*D_A^2Q(R^*D_A^2R)^{-1}J^*.
\end{equation}
To see this, f{}irst observe that the identity in (\ref{equal}) and
the def{}inition of $J$ imply that $Q^*D_A^2Q=R^*D_A^2R+J^*J$. With
this equality we obtain that
\begin{eqnarray*}
J^*\De_\Om
&=&J^*(I+J(R^*D_A^2R)^{-1}J^*)
=(I+J^*J(R^*D_A^2R)^{-1})J^*\\
&=&(R^*D_A^2R+Q^*D_A^2Q-R^*D_A^2R)(R^*D_A^2R)^{-1}J^*
=Q^*D_AQ(R^*D_AR)^{-1}J^*.
\end{eqnarray*}
So (\ref{K*DeOm}) holds.

{}From (\ref{K*DeOm}) it follows that
$J^*\De_\Om^{-1}=R^*D_A^2R(Q^*D_A^2Q)^{-1}J^*$ and
\[
(I-J(Q^*D_A^2Q)^{-1}J^*)\De_\Om
=\De_\Om-J(Q^*D_AQ)^{-1}J^*\De_\Om
=\De_\Om-J(R^*D_A^2R)^{-1}J^*=I.
\]
Since $\De_\Om$ is self adjoint, we see that the inverse of $\De_\Om$ is
given by
\begin{equation}\label{DeOminv}
\De_\Om^{-1}=I-J(Q^*D_A^2Q)^{-1}J^*.
\end{equation}
Moreover, we have $\Pi_{\sD_{T'}}J=D_{T'}AR$ and
$L_{D_AR}=(R^*D_A^2R)^{-1}R^*D_A$. Hence
\begin{eqnarray*}
L_{D_AR}^*J^*\De_\Om^{-1}JL_{D_AR}
&=&D_AR(R^*D_A^2R)^{-1}R^*D_A^2R(Q^*D_A^2Q)^{-1}J^*J
(R^*D_A^2R)^{-1}R^*D_A\\
&=&D_AR(Q^*D_A^2Q)^{-1}J^*J(R^*D_A^2R)^{-1}R^*D_A\\
&=&D_AR(Q^*D_A^2Q)^{-1}(Q^*D_A^2Q-R^*D_A^2R)(R^*D_A^2R)^{-1}R^*D_A\\
&=&D_AR((R^*D_AR)^{-1}-(Q^*D_A^2Q)^{-1})R^*D_A\\
&=&P_{\im D_AR}-D_AR(Q^*D_AQ)^{-1}R^*D_A.
\end{eqnarray*}

For the last equality note that, since $D_AR$ is left invertible, the
projection on $\im D_AR$ is given by
$P_{\im D_AR}=D_AR(R^*D_A^2R)^{-1}R^*D_A=D_ARL_{D_AR}$.
Using the formula for $\De_\Om^{-1}$ in (\ref{DeOminv}) and the fact
that $P_{\im D_AR}+P_{\kr R^*D_A}=I_{\sD_A}$ we obtain that
\begin{eqnarray*}
Y^*Y
&=&\mat{cc}{\Pi_{\sD_{T'}}\De_\Om^{-\half}&0\\
-L_{D_AR}^*J^*\De_\Om^{-\half}&-\Pi_{\kr R^*D_A}^*}
\mat{cc}{\De_\Om^{-\half}\Pi_{\sD_{T'}}^*&-\De_\Om^{-\half}JL_{D_AR}\\
0&-\Pi_{\kr R^*D_A}}\\
&=&\mat{cc}{\Pi_{\sD_{T'}}\De_\Om^{-1}\Pi_{\sD_{T'}^*}
&-\Pi_{\sD_{T'}}\De_\Om^{-1}JL_{D_AR}\\
-L_{D_AR}^*J^*\De_\Om^{-1}\Pi_{\sD_{T'}}^* &P_{\kr
R^*D_A}+L_{D_AR}^*J^*\De_\Om^{-1}JL_{D_AR}}\\
&=&\mat{cc}{\Pi_{\sD_{T'}}(I-J(Q^*D_A^2Q)^{-1}J^*)\Pi_{\sD_{T'}}^*
&-\Pi_{\sD_{T'}}J(Q^*D_A^2Q)^{-1}R^*D_A\\
-D_AR(Q^*D_A^2Q)^{-1}J^*\Pi_{\sD_{T'}} &I_{\sD_A}
-D_AR(Q^*D_A^2Q)^{-1}R^*D_A}\\
&=&\mat{cc}{I_{\sD_{T'}}&0\\0&I_{\sD_A}}-
\mat{c}{D_{T'}AR\\D_AR}(Q^*D_A^2Q)^{-1}\mat{cc}{R^*A^*D_{T'}&R^*D_A}
=D_{\om^*}^2.
\end{eqnarray*}
To see the last identity use the f{}irst formula for $\om\Pi_\sF$ in the
computation (\ref{omtrans}) and observe that
$L_{D_AQ}L_{D_AQ}^*=(Q^*D_A^2Q)^{-1}$.

Next we show that $\kr Y^*=\{0\}$.
Since $Y^*|\kr R^*D_A=-\Pi_{\kr R^*D_A}^*$ and
$Y^*(\sD_\circ\oplus\sD_{T'})\subset\sD_\circ\oplus\im D_AR
=\sD_\circ\oplus(\sD_A\ominus\kr R^*D_A)$, it suff{}ices to show that
$\kr (Y^*|\sD_\circ\oplus\sD_{T'})=\{0\}$. To see that this is the case
observe that $R^*D_AL_{D_AR}^*=I$. So $\kr L_{D_AR}^*=\{0\}$. Moreover,
we have $\kr D_\circ=\{0\}$, when $D_\circ$ is seen as an operator on
$\sD_\circ$. Therefore
\begin{eqnarray*}
\kr Y^*|\sD_\circ\oplus\sD_{T'}
&=&\kr\mat{c}{\Pi_{\sD_{T'}}\De_\Om^{-\half}\\
L_{D_AR}^*J^*\De_\Om^{-\half}}\\
&=&\kr\mat{cc}{I_{\sD_{T'}}&0\\0&L_{D_AR}^*}
\mat{c}{\Pi_{\sD_{T'}}\\J^*}\De_\Om^{-\half}\\
&=&\De_\Om^\half\kr\mat{c}{\Pi_{\sD_{T'}}\\J^*}
=\De_\Om^\half\kr\mat{cc}{0&I_{\sD_{T'}}\\D_\circ&R^*A^*D_{T'}}=\{0\}.
\end{eqnarray*}

\noindent{\bf Part 3.}
Let $N$ be an operator from $\sH_0$ to $\sH$ that has a left inverse.
In particular, this holds for $Q$ and $R$. Then $D_A^{-1}\Pi_{\kr N^*}$
has a left inverse, and thus
$\De_N=\Pi_{\kr N^*}D_A^{-2}\Pi_{\kr N^*}^*$ is positive def{}inite on
$\kr N^*$. Since $D_A^{-1}(D_AN)=N$, we see that
$D_A^{-1}\im D_AN=\im N$. This implies that
$\im D_A^{-1}\Pi_{\kr N^*}^*=D_A^{-1}\kr N^*=\kr N^*D_A$.
Then we obtain in the same way as for $D_AR$ in Part 2 that
\[
P_{\kr N^*D_A}=P_{\im D_A^{-1}\Pi_{\kr N^*}^*}
=D_A^{-1}\Pi_{\kr N^*}^*\De_N^{-1}\Pi_{\kr N^*}D_A^{-1}.
\]
Filling in $Q$ and $R$ for $N$ gives the desired results.\bigskip

\noindent{\bf Part 4.}
The identities $Y^*Y=D_{\om^*}$ and $\kr Y^*=\{0\}$ derived in Part 2
show that there exists a unitary operator $\phi_\om$ mapping
$\sD_{\om^*}$ onto $\sD_\circ\oplus\sD_{T'}\oplus\kr R^*D_A$ such that
$Y^*\phi_\om=D_{\om^*}$. Furthermore, from (\ref{projects}) we obtain
that there exist unitary operators $\phi_Q$ and $\phi_R$ mapping
$\kr Q^*$ and $\kr R^*$ onto $\kr Q^*D_A$ and $\kr R^*D_A$,
respectively, such that
$\phi_Q\De_Q^{-\half}\Pi_{\ker Q^*}D_A^{-1}=\Pi_{\kr Q^*D_A}$ and
$\phi_R\De_R^{-\half}\Pi_{\ker R^*}D_A^{-1}=\Pi_{\kr R^*D_A}$. Now put
\begin{equation*}\label{phi*}
\phi_*=\mat{cc}{I_{\sD_\circ\oplus\sD_{T'}}&0\\0&\phi_R^*}
\mat{c}{\Pi_{\sD_\circ\oplus\sD_{T'}}\phi_\om
\\\Pi_{\kr R^*D_A}\phi_\om}:\sD_{\om^*}\to
\mat{c}{\sD_\circ\oplus\sD_{T'}\\\kr R^*}.
\end{equation*}
Then $\phi_*$ is unitary and we have
\begin{eqnarray*}
\mat{c}{X_5\\D_AX_2}\phi_*
&=&\mat{cc}{\Pi_{\sD_{T'}}\De_\Om^{-\half}&0\\
-L_{D_AR}^*J^*\De_\Om^{-\half}
&-D_A^{-1}\Pi_{\kr R^*}^*\De_R^{-\half}}
\mat{cc}{I_{\sD_\circ\oplus\sD_{T'}}&0\\0&\phi_R^*}\times\\[.1cm]
&&\quad\quad\times\mat{c}{\Pi_{\sD_\circ\oplus\sD_{T'}}\phi_\om
\\\Pi_{\kr R^*D_A}\phi_\om}
=Y^*\phi_\om=D_{\om^*}.
\end{eqnarray*}
{}From the def{}inition of $\phi_Q$ we obtain that
\[
\phi_QX_3D_A^{-1}
=\phi_Q\De_Q^{-\half}\Pi_{\ker Q^*}D_A^{-1}=\Pi_{\kr Q^*D_A}.
\]
Finally, the fact that $\phi_Q$ and $\phi_*$ are unitary implies that
an operator-valued function $\Va$ is from the Schur class
$\eS(\kr Q^*,\sD_\circ\oplus\sD_{T'}\oplus\kr R^*D_A)$ if and only if
\begin{equation}\label{Vab}
\Va(\la)=\phi_*\Vb(\la)\phi_Q\quad\quad(\la\in\BD)
\end{equation}
for some Schur class function $\Vb$ from $\eS(\kr
Q^*D_A,\sD_{\om^*})$, and $\Va$ and $\Vb$ in (\ref{Vab}) def{}ine each
other uniquely. Therefore we obtain with the statement at the
beginning of the proof that $Z$ from $\eS(\sD_A,\sD_{T'}\oplus\sD_A)$
has $Z(\la)|\sF=\om$ for each $\la\in\BD$ if and only if there
exists a Schur class function $\Va$ from $\eS(\kr
Q^*,\sD_\circ\oplus\sD_{T'}\oplus\kr R^*D_A)$ such that
\begin{eqnarray*}
Z(\la)
&=&\om\Pi_\sF+D_{\om^*}\phi_*^*\Va(\la)\phi_Q^*\Pi_{\kr Q^*D_A}\\
&=&\mat{c}{X_4D_A^{-1}\\D_AX_1D_A^{-1}}+\mat{c}{X_5\\D_AX_2}
\Va(\la)X_3D_A^{-1}\quad\quad(\la\in\BD),
\end{eqnarray*}
and $\Va$ and $Z$ def{}ine each other uniquely.
\epr

\noindent{\bf Proof of Theorem \ref{mainth}}
Let $B$ be an operator from $\sH$ into $H^2(\sD_{T'})$. Assume that
$B$ is given by (\ref{solsA}) with $\tilGa$ the operator def{}ined by
(\ref{solsB}), where $Z$ is a function from the Schur class
$\eS(\sD_A,\sD_{T'}\oplus\sD_A)$ with $Z(\la)|\sF=\om$ for each
$\la\in\BD$. Let $\Va$ be the Schur class function from
$\eS(\kr Q^*,\sD_\circ\oplus\sD_{T'}\oplus\kr R^*D_A)$ such that $Z$ is
given by (\ref{ZV}). Put
\begin{equation}\label{XtilX}
\tilX_1=D_AX_1D_A^{-1},\quad
\tilX_2=D_AX_2,\quad
\tilX_3=X_3D_A^{-1},\quad
\tilX_4=X_4D_A^{-1}\ands
\tilX_5=X_5.
\end{equation}
Then
\begin{equation}\label{tilXrewrite}
\begin{array}{rcl}
\mat{c|c}{\tilX_3&0\\\hline\tilX_4&\tilX_5\\\tilX_1&\tilX_2}
&=&\mat{c|c}{X_3D_A^{-1}&0\\\hline X_4D_A^{-1}&X_5\\
D_AX_1D_A^{-1}&D_AX_2}\\[.2cm]
&=&\mat{cc}{\phi_Q^*&0\\0&I_{\sD_{T'}\oplus\sD_A}}
\mat{cc}{\Pi_{\sD_A\ominus\sF}&0\\\om\Pi_\sF&D_{\om^*}}
\mat{cc}{I_{\sD_A}&0\\0&\phi_*^*}.
\end{array}
\end{equation}
Here $\phi_Q$ and $\phi_*$ are the unitary operators constructed in
Part 4 of the proof of Proposition \ref{pr:BTR}. Since $\om$ is a
contraction, we see that the operator
\begin{equation}\label{tilXuneq}
\mat{cc}{\Pi_{\sD_A\ominus\sF}&0\\\om\Pi_\sF&D_{\om^*}}
:\mat{c}{\sD_A\\\sD_{\om^*}}\to\mat{c}{\kr Q^*D_A\\\sD_{T'}\oplus\sD_A}
\end{equation}
is a contraction. After rearranging rows in (\ref{tilXrewrite}) we
obtain, from the fact that $\psi_Q$ and $\psi_*$ are unitary, that
$\tilX$ in (\ref{tilX}) is a contraction.
Let $\Pb_{11}$, $\Pb_{12}$, $\Pb_{21}$ and $\Pb_{22}$
be the functions def{}ined in (\ref{Pb}) and Let $\Pa_{11}$,
$\Pa_{12}$, $\Pa_{21}$ and $\Pa_{22}$ be the functions def{}ined in
(\ref{sols3}). Using that $(I-\la\tilX_1)^{-1}=D_A(I-\la
X_1)^{-1}D_A^{-1}$ we obtain from the relations in (\ref{XtilX}) that
for each $\la\in\BD$
\begin{equation}\label{PbPa}
\begin{array}{l}
\Pb_{11}(\la)=\Pa_{11}(\la),\quad
\Pb_{12}(\la)=\Pa_{12}(\la)D_A^{-1},\\[.2cm]
\Pb_{21}(\la)=\Pa_{21}(\la)\ands
\Pb_{22}(\la)=\Pa_{22}(\la)D_A^{-1}.
\end{array}
\end{equation}
Applying Proposition \ref{pr:RED} we see that $\Pa_{11}$ and
$\Pa_{21}$ are Schur class functions from
$\eS(\sD_\circ\oplus\sD_{T'}\oplus\kr R^*,\kr Q^*)$ and
$\eS(\sD_\circ\oplus\sD_{T'}\oplus\kr R^*,\sD_{T'})$, respectively,
and that $\Pa_{12}$ and $\Pa_{22}$ are functions from $\eH^2(\sH,\kr Q^*)$
and $\eH^2(\sH,\sD_{T'})$, respectively. Moreover, we have for each
$h\in\sH$ and each $\la\in\BD$ that
\begin{eqnarray*}
(\tilGa D_Ah)(\la)&=&\Pi_{\sD_{T'}} Z(\la)(I-\la\Pi_{\sD_A}Z(\la))D_Ah\\
&=&\Pb_{22}(\la)D_Ah+\Pb_{21}(\la)\Va(\la)(I-\Pb_{11}(\la)\Va(\la))^{-1}
\Pb_{12}(\la)D_Ah\\
&=&\Pa_{22}(\la)h+\Pa_{21}(\la)\Va(\la)(I-\Pa_{11}(\la)\Va(\la))^{-1}
\Pa_{12}(\la)h=(\Ga h)(\la),
\end{eqnarray*}
where $\Ga$ is the operator def{}ined in (\ref{sols2}) with $\Va$ the
Schur class function determined by (\ref{ZV}).
\epr

%%%%%%%%%%%%%%%%%%%%%%%%%%%%%%%%%%%%%%%%%%%%%%%%%%%%%%%%%%%%%%%%%%%%%%%%
\subsection{Corollaries}

We conclude this section with two corollaries. The f{}irst
specif{}ies Theorem \ref{mainth} for the classical commutant lifting
setting.

%%%%%%%%%%%%%%%%%%%%%%%%%%%%%%%%%%%
\begin{corollary}\label{cor:CCL}
Let $\LDS$ be a lifting data set where $U'$ is the
Sz.-Nagy-Sch\"affer isometric lifting of $T'$. Assume that
$Q$ is an isometry on $\sH$, $R=I_\sH$ and $A$ is a strict contraction.
Then $\kr R^*=\{0\}$, $\sD_\circ=\{0\}$ and the operator-valued
functions $\Pa_{11}$, $\Pa_{12}$, $\Pa_{21}$ and $\Pa_{22}$ in
\tu{(\ref{sols3})} are functions from $\eS(\sD_{T'},\kr Q^*)$,
$\eH^\infty(\sH,\kr Q^*)$, $\eS(\sD_{T'},\sD_{T'})$ and
$\eH^\infty(\sH,\sD_{T'})$, respectively, and they are given by
\begin{equation}\label{PaCCL}
\begin{array}{rcl}
\Pa_{11}(\la)&=&-\la\De_Q^{-\half}\Pi_{\kr Q^*}
(I-\la T_A)^{-1}D_A^{-2}A^*D_{T'}\De_\Om^{-\half},\\[.2cm]
\Pa_{12}(\la)&=&\De_Q^{-\half}\Pi_{\kr Q^*}(I-\la T_A)^{-1},\\[.2cm]
\Pa_{21}(\la)&=&\De_\Om^\half-D_{T'}A(I-\la T_A)^{-1}
D_A^{-2}A^*D_{T'}\De_\Om^{-\half},\\[.2cm]
\Pa_{22}(\la)&=&D_{T'}A(I-\la T_A)^{-1}T_A,
\end{array}
\quad\quad\quad(\la\in\BD)
\end{equation}
where
\begin{equation}\label{TA}
T_A=D_{AQ}^{-2}Q^*D_A^2\ons\sH,
\end{equation}
\begin{equation}\label{CCL:Deltas}
\De_Q=\Pi_{\kr Q^*}D_A^{-1}\Pi_{\kr Q^*}^*\ons\kr Q^* \ands
\De_\Om=I+D_{T'}AD_A^{-2}A^*D_{T'}\ons\sD_{T'}.
\end{equation}
\end{corollary}

\bpr
Since $R=I_\sH$ and $R^*R=I_\sH=Q^*Q$, we obtain that $\kr
R^*=\{0\}$ and $D_\circ=0$ on $\sD_\circ=\{0\}$. Clearly $R$ is left
invertible and, by assumption, $A$ is a strict contraction.
Moreover, we have that $R-\la Q=I_\sH-\la Q$ is (left) invertible
for each $\la\in\BD$, because $Q$ is contractive. Hence the result
of Theorem \ref{mainth}, and the remark in the third paragraph
underneath Theorem \ref{mainth} hold for this lifting data set. In
particular $\Pa_{11}$, $\Pa_{12}$, $\Pa_{21}$ and $\Pa_{22}$ in
\tu{(\ref{sols3})} are functions from $\eS(\sD_{T'},\kr Q^*)$,
$\eH^\infty(\sH,\kr Q^*)$, $\eS(\sD_{T'},\sD_{T'})$ and
$\eH^\infty(\sH,\sD_{T'})$, respectively.

Observe that $J$ in (\ref{Deltas})
is given by $J=D_{T'}A$ and we have
\[
Q^*D_A^2Q=Q^*Q-Q^*A^*AQ=I-Q^*A^*AQ=D_{AQ}^2.
\]
Therefore we have that $X_1=T_A$ with $X_1$ from (\ref{X1X5}) and $T_A$
as in (\ref{TA}). Note that $\De_Q$ and $\De_\Om$ given in
(\ref{Deltas}) reduce to the formulas in
(\ref{CCL:Deltas}) So we obtain that the operators $X_1$, $X_2$, $X_3$,
$X_4$ and $X_5$ in (\ref{X1X5}), under the present assumptions, are
\[
X_1=T_A,\quad\!
X_2=-D_A^{-2}A^*D_{T'}\De_\Om^{-\half},\quad\!
X_3=\De_Q^{-\half}\Pi_{\kr Q^*},\quad\!
X_4=D_{T'}A T_A\!\ands\!
X_5=\De_\Om^{-\half}.
\]
This immediately shows that the formulas for $\Pa_{11}$, $\Pa_{12}$ and
$\Pa_{22}$ in (\ref{sols3}) are given by (\ref{PaCCL})
in the classical commutant lifting setting. Furthermore, we have for
each $\la\in\BD$ that
\begin{eqnarray*}
\Pa_{21}(\la)
&=&X_5+\la X_4(I-\la X_1)^{-1}X_2
=\De_\Om^{-\half}
-\la D_{T'}AT_A(I-\la T_A)^{-1}D_A^{-1}A^*D_{T'}\De_\Om^{-\half}\\
&=&(I+D_{T'}AD_A^{-1}A^*D_{T'})\De_\Om^{-\half}
-D_{T'}A(I-\la T_A)^{-1}D_A^{-1}A^*D_{T'}\De_\Om^{-\half}\\
&=&\De_\Om\De_\Om^{-\half}
-D_{T'}A(I-\la T_A)^{-1}D_A^{-1}A^*D_{T'}\De_\Om^{-\half}\\
&=&\De_\Om^\half-D_{T'}A(I-\la T_A)^{-1}D_A^{-1}A^*D_{T'}
\De_\Om^{-\half}.
\end{eqnarray*}
So $\Pa_{21}$ in (\ref{sols3}) also reduces to its formula in
(\ref{PaCCL}).
\epr

%%%%%%%%%%%%%%%%%%%%%%%%%%%%%%%%%%%
\begin{corollary}\label{cor:multop}
Let $\Om=\LDS$ be a lifting data set with $U'$ the Sz.-Nagy-Sch\"affer
isometric lifting of $T'$.  Assume that $A$ is a strict contraction and
$R$ has a left inverse. Let $\Pa_{11}$, $\Pa_{12}$, $\Pa_{21}$ and
$\Pa_{22}$ be the functions def{}ined by \tu{(\ref{sols3})} and
\tu{(\ref{X1X5})}. Put
\begin{equation}\label{multop}
M=\mat{cc}{0&A\\M_{\Pa_{11}}&\Ga_{\Pa_{12}}\\
M_{\Pa_{21}}&\Ga_{\Pa_{22}}}
:\mat{c}{H^2(\sD_\circ\oplus\sD_{T'}\oplus\kr R^*)\\\sH}\to
\mat{c}{\sH'\\H^2(\kr Q^*)\\H^2(\sD_{T'})}.
\end{equation}
Then $M$ is a contraction, and $M$ is an isometry if $R^*R=Q^*Q$ and
$X_1$ is pointwise stable.
\end{corollary}

\epr
Let $X_1$, $X_2$, $X_3$, $X_4$ and $X_5$ be the operators in
(\ref{X1X5}) and $\tilX_1$, $\tilX_2$, $\tilX_3$,  $\tilX_4$ and
$\tilX_5$ the operators given by (\ref{XtilX}). Def{}ine
$\Pb_{11}$, $\Pb_{12}$, $\Pb_{21}$ and $\Pb_{22}$ by (\ref{Pb}).
We saw in the proof of Theorem \ref{mainth} that the results from
Proposition \ref{pr:RED} can be applied for this choice of $\tilX_1$,
$\tilX_2$, $\tilX_3$,  $\tilX_4$ and  $\tilX_5$. Then the operator
$\tilM$ in (\ref{tilmultop}) is a contraction which is unitary in case
$\tilX$ in (\ref{tilX}) is unitary and $\tilX_1$ is pointwise stable.
Since $X_1$ and $\tilX_1$ are similar ($D_AX_1=\tilX_1D_A$ and $D_A>0$),
we see that $\tilX_1$ is pointwise stable if and only if $X_1$ is
pointwise stable. Moreover, the computation in (\ref{tilXrewrite}) and
the remark below (\ref{tilXuneq}) show that $\tilX$ is unitary if and
only if the operator in (\ref{tilXuneq}) is unitary. This is precisely
the case if and only if $\om$ is an isometry. In other words, $\tilX$
in (\ref{tilX}) is unitary if and only if $R^*R=Q^*Q$. We conclude that
the operator $\tilM$ in (\ref{tilmultop}) is a contraction which is
unitary if $X_1$ is pointwise stable and $R^*R=Q^*Q$.

In terms of operators the identities in (\ref{PbPa}) become
\begin{equation}\label{opidens}
M_{\Pb_{11}}=M_{\Pa_{11}},\quad
\Ga_{\Pb_{12}}=\Ga_{\Pa_{11}}D_A^{-1},\quad
M_{\Pb_{21}}=M_{\Pa_{21}}\ands
\Ga_{\Pb_{11}}=\Ga_{\Pa_{11}}D_A^{-1}.
\end{equation}
Next def{}ine a contraction $\ov{A}$ by
\[
\ov{A}=\mat{cc}{0&A}:\mat{c}{H^2(\sW)\\\sH}\to\sH',
\]
where $\sW=\sD_\circ\oplus\sD_{T'}\oplus\kr R^*$.
The defect operator of $\ov{A}$ is given by
\[
D_{\ov{A}}=\mat{cc}{I_{H^2(\sW)}&0\\0&D_A}
:\mat{c}{H^2(\sW)\\\sH}\to\mat{c}{H^2(\sW)\\\sD_A}.
\]
Then the identities in (\ref{opidens}) show that
\begin{equation*}\label{CtilC}
M=\mat{cc}{0&A\\M_{\Pa_{11}}&\Ga_{\Pa_{12}}\\
M_{\Pa_{21}}&\Ga_{\Pa_{22}}}
=\mat{cc}{0&A\\M_{\Pb_{11}}&\Ga_{\Pb_{12}}D_A\\
M_{\Pb_{21}}&\Ga_{\Pb_{22}}D_A} =\mat{c}{\ov{A}\\\tilM D_{\ov{A}}}
=\mat{cc}{I_{\sH'}&0\\0&\tilM}\mat{c}{\ov{A}\\D_{\ov{A}}}.
\end{equation*}
It is well known that the operator $\mat{cc}{\ov{A}^*&D_{\ov{A}}}^*$
is an isometry. Therefore we obtain that $M$ is a contraction which is
isometric in case $R^*R=Q^*Q$ and $X_1$ is pointwise stable.
\epr

%%%%%%%%%%%%%%%%%%%%%%%%%%%%%%%%%%%%%%%%%%%%%%%%%%%%%%%%%%%%%%%%%%%%%%%%
%%%%%%%%%%%%%%%%%%%%%%%%%%%%%%%%%%%%%%%%%%%%%%%%%%%%%%%%%%%%%%%%%%%%%%%%
\section{The relaxed Nehari extension problem}
\label{sec:Nehari}

In this section we apply Theorem \ref{mainth} to obtain a description
of all $N$-complementary sequences for the relaxed Nehari extension
problem formulated in the last but one paragraph of the introduction.

%%%%%%%%%%%%%%%%%%%%%%%%%%%%%%%%%%%%%%%%%%%%%%%%%%%%%%%%%%%%%%%%%%%%%%%%
\subsection{The relaxed Nehari problem in the relaxed commutant
lifting setting}

Throughout this subsection $N$ is a positive integer and
$F_{-1},F_{-2},\ldots$ is a sequence of operators from $\sU$ to
$\sY$ satisfying $\sum_{n=1}^\infty\|F_{-n}u\|^2<\infty$ for each
$u\in\sU$. We refer to the operator from $\sU^N$ into
$\ell^2_-(\sY)$ given by
\begin{equation}
\label{NHankel2}\mat{cccc}
{\vdots&\vdots&&\vdots\\
F_{-3}&F_{-4}&\cdots&F_{-(N+2)}\\
F_{-2}&F_{-3}&\cdots&F_{-(N+1)}\\
F_{-1}&F_{-2}&\cdots&F_{-N}}: \sU^N\to \ell^2_-(\sY)
\end{equation}
as the \emph{$N$-truncated Hankel operator} def{}ined by the sequence
$(F_{-n-1})_{n\in\BN}$. Our f{}irst remark is that an operator $A$
from $\sU^N$ into $\ell^2_-(\sY)$ is an $N$-truncated Hankel
operator if and only if $A$ satisf{}ies the intertwining relation
$T'AR=AQ$, where $T'$ is the shift operator on $\ell^2_-(\sY)$ given
by
\begin{equation}
\label{defTpr}
T'=\mat{ccccc} {\ddots&&\vdots&\vdots&\vdots\\
\ddots&\ddots&\vdots&\vdots&\vdots\\\ddots&\ddots&0&0&0\\
\cdots&0&I_\sY&0&0\\\cdots&0&0&I_\sY&0},
\end{equation}
and $R$ and $Q$ are the operators from $\sU^{N-1}$ to $\sU^N$
def{}ined by
\begin{equation}\label{RandQ}
R=\mat{c}{I_{\sU^{N-1}}\\0}:\sU^{N-1}\to\mat{c}{\sU^{N-1}\\\sU}\ands
Q=\mat{c}{0\\I_{\sU^{N-1}}}:\sU^{N-1}\to\mat{c}{\sU\\\sU^{N-1}}.
\end{equation}
Indeed, if $A$ is an operator from $\sU^N$ into $\ell^2_-(\sY)$,
then $A$ admits an operator matrix decomposition of the form
\[
A=\mat{cccc}
{\vdots&\vdots&&\vdots\\
A_{31}&A_{32}&\cdots&A_{3N}\\
A_{21}&A_{22}&\cdots&A_{2N}\\
A_{11}&A_{12}&\cdots&A_{1N}},
\]
and hence
\[
AR=\mat{cccc}
{\vdots&\vdots&&\vdots\\
A_{31}&A_{32}&\cdots&A_{3\, N-1}\\
A_{21}&A_{22}&\cdots&A_{2\, N-1}\\
A_{11}&A_{12}&\cdots&A_{1\, N-1}}\ands AQ=\mat{cccc}
{\vdots&\vdots&&\vdots\\
A_{32}&A_{33}&\cdots&A_{3N}\\
A_{22}&A_{23}&\cdots&A_{2N}\\
A_{12}&A_{13}&\cdots&A_{1N}}.
\]
It follows that $T'AR=AQ$ is equivalent to $A_{kj}=A_{k-1\,j+1}$ for
appropriate indices $k$ and $j$, that is, $T'AR=AQ$ is equivalent to
$A$ being an $N$-truncated Hankel operator.

We shall also need the bilateral forward shift $V'$ on $\ell^2(\sY)$
which is given by
\begin{equation}\label{defVpr}
V'=\mat{ccccccccc}
{\ddots&\ddots&\ddots&\ddots\ddots\\&0&I&0&0&0\\&&0&I&\fbox{0}&0&0\\
&&&0&I&0&0&0\\&&&&\ddots&\ddots&\ddots&\ddots&\ddots}
\ons\ell^2(\sY).
\end{equation}
Let $\ell^2_+(\sY)$ be the subspace of $\ell^2(\sY)$ consisting of the
sequences that have the zero vector on all entries with strictly
negative index. Identifying $\ell^2_-(\sY)$ with the subspace of
$\ell^2(\sY)$ consisting of all sequences that have the zero vector
on all entries with a non-negative index, we see that
$\ell^2(\sY)=\ell^2_-(\sY)\oplus \ell^2_+(\sY)$, and $V'$ partitions
as
\begin{equation}
\label{partVpr} V'= \mat{cc}{T'&0\\X'&S'}
:\mat{c}{\ell^2_-(\sY)\\\ell^2_+(\sY)}\to \mat{c}{\ell^2_-(\sY)\\
\ell^2_+(\sY)}.
\end{equation}
Here $T'$ is given by (\ref{defTpr}), the operator $S'$ is the
unilateral forward shift on $\ell^2_+(\sY)$, and $X'$ is the
operator from $\ell^2_-(\sY)$ into $\ell^2_+(\sY)$ given by
\begin{equation}\label{defXpr}
X'(\ldots, y_{-3}, y_{-2},y_{-1})=(y_{-1}, 0, 0,\ldots)\quad\quad
((\ldots, y_{-3}, y_{-2},y_{-1})\in\ell^2_-(\sY)).
\end{equation}
Since $V'$ is unitary, (\ref{partVpr}) shows that $V'$ is an
isometric lifting of $T'$. As is easily seen this lifting is also
minimal.

We are now ready to state the main result of this subsection.

%%%%%%%%%%%%%%%%%%%%%%%%%%%%%%%%%%%
\begin{proposition}\label{pr:Neh}
Let $N$ be a positive integer and $F_{-1},F_{-2},\ldots$  a sequence
of operators from $\sU$ to $\sY$ satisfying
$\sum_{n=1}^\infty\|F_{-n}u\|^2<\infty$ for each $u\in\sU$. Assume
 that the $N$-truncated Hankel operator $A$
def{}ined by  $(F_{-n-1})_{n\in\BN}$ is a contraction. Then
$\tilde\Om=\{A,T',V',R,Q\}$, with $T',V',R$ and $Q$ def{}ined by
\tu{(\ref{defTpr})}, \tu{(\ref{defVpr})} and \tu{(\ref{RandQ})}, is a lifting
data set, and there exists a contractive interpolant for $\tilde\Om$
if and only if there exists an $N$-complementary sequence associated
with $(F_{-n-1})_{n\in\BN}$. More precisely, an operator $\tilB$ from
$\sU^N$ to $\ell^2(\sY)$ is a contractive interpolant for $\tilde\Om$
if and  only if $\tilB$ is of the form
\begin{equation}\label{defB}
\tilB=\mat{cccc}
{\vdots&\vdots&&\vdots\\
F_{-2}&F_{-3}&\cdots&F_{-(N+1)}\\
F_{-1}&F_{-2}&\cdots&F_{-N}\\
\fbox{$H_0$}
&F_{-1}&\cdots&F_{-(N-1)}\\
H_1&H_0&\cdots &\vdots\\
\vdots&\ddots&\ddots&F_{-1}\\
\vdots&&\ddots&H_0\\
\vdots&&&\vdots\\}:\sU^N\to\ell^2(\sY)
\end{equation}
with  $(H_n)_{n\in\BN}$ an $N$-complementary sequence associated
with $(F_{-n-1})_{n\in\BN}$.
\end{proposition}

\bpr
We already know that $T'AR=AQ$ and that $V'$ is a minimal
isometric lifting of $T'$. Since $R$ and $Q$ are both isometries, we
have $R^*R=Q^*Q$. So the constraints (\ref{lds}) are fulf{}illed,
and hence $\tilde\Om$ is a lifting data set.

Now let $\tilB$ be an operator from $\sU^N$ into $\ell^2(\sY)$. Using
similar arguments as in the f{}irst paragraph of this subsection we
obtain that $\tilB$ satisf{}ies $V'\tilB R=\tilB Q$ if and only if
\begin{equation}\label{defB2}
\tilB=\mat{cccc}
{\vdots&\vdots&&\vdots\\
H_{-2}&H_{-3}&\cdots&H_{-(N+1)}\\
H_{-1}&H_{-2}&\cdots&H_{-N}\\
\fbox{$H_0$}
&H_{-1}&\cdots&H_{-(N-1)}\\
H_1&H_0&\cdots&H_{-(N-2)}\\
H_2&H_1&\cdots&H_{-(N-3)}\\
\vdots&\vdots&&\vdots\\}:\sU^N\to\ell^2(\sY)
\end{equation}
for some sequence of operators
$(H_n)_{n\in\BZ}$
%$\dots, H_{-2},H_{-1}, H_0, H_1, H_2, \ldots$
from $\sU$ into
$\sY$ satisfying $\sum_{n=-\infty}^\infty\|H_{n}u\|^2<\infty$ for
each $u\in\sU$. Moreover, if $\tilB$ is given by (\ref{defB2}), then
$\Pi_{\ell^2_-(\sY)}\tilB=A$ holds if and only if $H_k=F_k$ for
$k=-1,-2,-3,\ldots$. Henceforth, we obtain that $\tilB$ is a
contractive  interpolant for $\tilde\Om$ if and only if $\tilB$ is
given by (\ref{defB})  with $(H_n)_{n\in\BN}$ an $N$-complementary
sequence associated with $(F_{-n-1})_{n\in\BN}$.
\epr

In order to apply Theorem \ref{mainth} we need to establish a relation
between the contractive interpolants for the lifting data set
$\tilde\Om=\{A,T',V',R,Q\}$ in Proposition \ref{pr:Neh} and those for
$\Om=\LDS$, with $U'$ being the Sz.-Nagy-Sch\"affer isometric lifting
of $T'$. For that purpose we use the Fourier transform $\Four_\sY$
from $\ell^2_+(\sY)$ to the Hardy space $H^2(\sY)$, that is,
$\Four_\sY$ is the unitary operator mapping $\ell^2_+(\sY)$ onto
$H^2(\sY)$ def{}ined by
\begin{equation}\label{four}
(\Four_\sY(y_n)_{n\in\BN})(\la)=\sum_{n=0}^\infty\la^ny_n\quad\quad
((y_n)_{n\in\BN}\in\ell^2_+(\sY),\ \la\in\BD).
\end{equation}
Using $\Four_\sY$, let $\Psi$ be the unitary operator from
$\ell^2(\sY)=\ell^2_-(\sY)\oplus \ell^2_+(\sY)$ into
$\ell^2_-(\sY)\oplus H^2(\sY)$  given by
\begin{equation}\label{uneqlifts}
\Psi=\mat{cc}{I_{\displaystyle{\ell^2_-(\sY)}}&0\\0&\Four_\sY}
:\mat{c}{\ell^2_-(\sY)\\\ell^2_+(\sY)}\to \mat{c}{\ell^2_-(\sY)\\
H^2(\sY)}.
\end{equation}
Since $\sD_{T'}=\sY$ and the bilateral shift $S'$ on $\ell^2_+(\sY)$
and the bilateral shift $S_\sY$ on $H^2(\sY)$ are related via
$\Four_\sY S'=S_\sY\Four_\sY$, we see from (\ref{partVpr}) and
(\ref{szns}) that $\Psi$ intertwines the bilateral shift $V'$ with
the Sz.-Nagy-Sch\"affer isometric lifting $U'$ of $T'$, that is,
$\Psi V'=U'\Psi$. Moreover, we have that $\Psi f=f$ for each
sequence $f$ in $\ell^2_-(\sY)$. Using the properties of $\Psi$
found above, a straightforward computation shows that $B$ is a
contractive interpolant for $\Om$ if and only if $B=\Psi\tilB$ for
some contractive interpolant $\tilB$ for $\tilde\Om$.

Finally, observe that a sequence of operators $H_0,H_1,\ldots$ in
$\Ops(\sU,\sY)$ has the property that
$\sum_{n=0}^\infty\|H_nu\|^2<\infty$ for each $u\in\sU$ if and only if
it is the sequence of Taylor coeff{}icients at zero of a function in
$\eH^2(\sU,\sY)$. So, alternatively, we seek functions $H$ in
$\eH^2(\sU,\sY)$ such that $\tilB$ in (\ref{defB}) is a contraction,
where $H_n$ is the $n^\tu{th}$ Taylor coeff{}icient of $H$ at zero.

%%%%%%%%%%%%%%%%%%%%%%%%%%%%%%%%%%%%%%%%%%%%%%%%%%%%%%%%%%%%%%%%%%%%%%%%
\subsection{The solution to the relaxed Nehari problem}

In the previous subsection we saw how the solution to the relaxed
commutant lifting problem can be applied to obtain all
$N$-complementary sequences for the relaxed Nehari extension problem.
Two important operators in the description of all contractive
interpolants in Theorem \ref{mainth} are the defect operators of
the contractions $T'$ and $A$. {}From the def{}inition of $T'$ in
(\ref{defTpr}) we immediately see that
\begin{equation}\label{defectT'}
\sD_{T'}=\sY\ands D_{T'}=\mat{cccc}{\cdots&0&0&I_\sY}
:\ell^2_-(\sY)\to\sY.
\end{equation}
The defect operator of the $N$-truncated Hankel operator, denoted by
$A$, is the positive square root of
\begin{equation}\label{defectA}
D_A^2=\mat{ccc}{\La_{1,1}&\cdots&\La_{1,N}\\\vdots&\ddots&\vdots\\
\La_{N,1}&\cdots&\La_{N,N}}\ons\sU^N,\mbox{ where }
\La_{i,j}=
\left\{
\begin{array}{c}
I-\sum_{n=i}^\infty F_{-n}^*F_{-n}\ \mbox{ if }i=j,\\[.2cm]
-\sum_{n=i}^\infty F_{-n}^*F_{-n+i-j}\mbox{ if }i\not=j.
\end{array}
\right.
\end{equation}

In what follows we assume the $N$-truncated Hankel operator to
be a strict contraction, or equivalently, we assume that the defect
operator $D_A$ is positive def{}inite. We also use the entries
$\Lat_{i,j}$ in the $N\times N$ operator matrix representation of
$D_A^{-2}$,
\begin{equation}\label{defectAinv}
D_A^{-2}=\mat{ccc}{\Lat_{1,1}&\cdots&\Lat_{1,N}\\
\vdots&\ddots&\vdots\\\Lat_{N,1}&\cdots&\Lat_{N,N}}\ons\sU^N.
\end{equation}
The fact that $D_A^2$ is positive def{}inite implies that $D_A^{-2}$
and $\Lat_{n,n}$ in (\ref{defectAinv}) are also positive def{}inite
for $n=1,\ldots,N$. In particular, this is true for $n=1$ and $n=N$.
Moreover, we have that the $N-1$ by $N-1$ left upper block matrix
operator of $D_A^2$ in (\ref{defectA}) is positive def{}inite.
Therefore there exists a unique solution
$\mat{ccc}{G_1&\cdots&G_{N-1}}$ to the equation
\begin{equation}\label{G1GN-1}
\mat{ccc}{\La_{1,1}&\cdots&\La_{1,N-1}\\\vdots&\ddots&\vdots\\
\La_{N-1,1}&\cdots&\La_{N-1,N-1}}\mat{c}{G_1^*\\\vdots\\G_{N-1}^*}
=\mat{c}{F_{-1}^*\\\vdots\\F_{-N+1}^*}.
\end{equation}

A description of all $N$-complementary sequences associated with
$(F_{-n-1})_{n\in\BN}$ under the assumption that the
$N$-truncated Hankel operator for $(F_{-n-1})_{n\in\BN}$ is a strict
contraction is given in the next theorem.

%%%%%%%%%%%%%%%%%%%%%%%%%%%%%%%%%%%
\begin{theorem}\label{th:Neh}
Let $N$ be a positive integer and $F_{-1},F_{-2},\ldots$  a sequence
of operators from $\sU$ to $\sY$ satisfying
$\sum_{n=1}^\infty\|F_{-n}u\|^2<\infty$ for each $u\in\sU$. Assume
that the $N$-truncated Hankel operator $A$ def{}ined by
$(F_{-n-1})_{n\in\BN}$ is a contraction. Let $H$ be a
function from $\eH^2(\sU,\sY)$. Then the Taylor coeff{}icients
$(H_n)_{n\in\BN}$ of $H$ at zero form an $N$-complementary sequence
associated with $(F_{-n-1})_{n\in\BN}$ if and only if there exists a
Schur class function $\Va$ from $\eS(\sU,\sY\oplus\sU)$ such that
\begin{equation}\label{NehPhi1}
H(\la)=\Pc_{22}(\la)+\Pc_{21}(\la)\Va(\la)
(I-\Pc_{11}(\la)\Va(\la))^{-1}\Pc_{12}(\la)\quad\quad(\la\in\BD).
\end{equation}
Moreover, for any function $V$ from $\eS(\sU,\sY\oplus\sU)$ the formula
\tu{(\ref{NehPhi1})} def{}ines a function $H$ from $\eH^2(\sU,\sY)$.
Here $\Pc_{11}$ and $\Pc_{21}$ are Schur class functions from
$\eS(\sY\oplus\sU,\sU)$ and $\eS(\sY\oplus\sU,\sY)$, respectively,
and $\Pc_{12}$ and $\Pc_{22}$ are functions from $\eH^\infty(\sU,\sU)$
and $\eH^\infty(\sU,\sY)$, respectively, and these functions are given
by
\begin{equation}\label{NehPhi2}
\begin{array}{rcl}
\Pc_{11}(\la)&=&-\la(\Lat_{1,1})^{-\half}E_N^*
(I_{\sU^N}-\la \TA)^{-1}\mat{cc}{G^*(I+FG^*)^{-\half}&C_2},\\[.2cm]
\Pc_{12}(\la)&=&(\Lat_{1,1})^{-\half}
-\la(\Lat_{1,1})^{-\half}E_N^*(I_{\sU^N}-\la \TA)^{-1}C_1,\\[.2cm]
\Pc_{21}(\la)&=&\mat{cc}{(I+FG^*)^\half&FC_2}
-F(I_{\sU^N}-\la \TA)^{-1}\mat{cc}{G^*(I+FG^*)^{-\half}&C_2},\\[.2cm]
\Pc_{22}(\la)&=&-F(I_{\sU^N}-\la \TA)^{-1}C_1,
\end{array}
\end{equation}
where $\TA$ on $\sU^N$, $E_N:\sU\to\sU^N$, $C_1:\sU\to\sU^N$,
$C_2:\sU\to\sU^N$, $F:\sU^N\to\sY$ and $G:\sU^N\to\sY$ are the
operators given by
\begin{equation}\label{NehCons}
\begin{array}{l}
\TA=\mat{ccccc}
{-\Lat_{2,1}(\Lat_{1,1})^{-1}&I&0&\cdots&0\\
-\Lat_{3,1}(\Lat_{1,1})^{-1}&0&I&\ddots&\vdots\\
\vdots&\vdots&\ddots&\ddots&0\\
-\Lat_{N,1}(\Lat_{1,1})^{-1}&0&&\ddots&I\\
0&0&0&\cdots&0},\quad
E_N=\mat{c}{I\\0\\\vdots\\\vdots\\0},\\[.2cm]
C_1=\mat{c}{\Lat_{2,1}\\\vdots\\\Lat_{N,1}\\0}
(\Lat_{1,1})^{-1},\quad
C_2=\mat{c}{\Lat_{1,N}\\\vdots\\\Lat_{N,N}}
(\Lat_{N,N})^{-\half},\\[.2cm]
F=\mat{ccc}{F_{-1}&\cdots&F_{-N}}\ands
G=\mat{cccc}{G_1&\cdots&G_{N-1}&0},
\end{array}
\end{equation}
with $\Lat_{m,n}$ for $m,n=1,\ldots,N$ as in \tu{(\ref{defectAinv})}, and
$\mat{ccc}{G_1&\cdots&G_{N-1}}$ the solution to the equation
\tu{(\ref{G1GN-1})}. Furthermore, we have $\spec(\TA)<1$, and the
operator $\hat M$ def{}ined by
\begin{equation}\label{hatM}
\hat M=\mat{cc}{0&\tilGa_-\\
M_{\Pc_{11}}&\Ga_{\Pc_{12}}\\
M_{\Pc_{21}}&\Ga_{\Pc_{22}}}
:\mat{c}{H^2(\sY\oplus\sU)\\\sU}\to
\mat{c}{\ell^2_-(\sY)\\H^2(\sU)\\H^2(\sY)},
\end{equation}
is isometric, where $\tilGa_-$ is given by
$\tilGa_-^*=\mat{ccc}{\cdots&F_{-2}^*&F_{-1}^*}:\ell^2_-(\sY)\to\sU$.
\end{theorem}

Observe that the state operator $\TA$ in (\ref{NehCons}) is close to
a companion operator. To be precise, let $E$ be the f{}lip over
operator on $\sU^N$ given by
\[
E(u_1,u_2,\ldots,u_N)=(u_N,u_{N-1},\ldots,u_1)\quad\quad((u_1,u_2,\ldots,u_N)\in\sU^N).
\]
Then $E$ is unitary and $E\TA E$ is precisely the second companion
operator corresponding to the operator-valued polynomial
$K(\la)=\la\sum_{k=0}^{N-1}\la^k\Lat_{N-k,1}$, see Chapter 14 in
\cite{LT85}. Note that the leading coeff{}icient $\Lat_{1,1}$ of $K$
is positive def{}inite.

As a computational remark, note that to obtain the solution
$\mat{ccc}{G_1&\cdots&G_{N-1}}$ to the equation (\ref{G1GN-1}) it
suff{}ices to compute the inverse of the operator $\Lat_{N,N}$
in (\ref{defectAinv}). Indeed, with a standard Schur complement type
of argument we see that the inverse of the $N-1$ by $N-1$ left upper
corner of $D_A^2$ in (\ref{defectA}) is given by
\begin{equation*}
\mat{ccc}{\Lat_{1,1}&\cdots&\Lat_{1,N-1}\\
\vdots&\ddots&\vdots\\\Lat_{N-1,1}&\cdots&\Lat_{N-1,N-1}}
-\mat{c}{\Lat_{1,N}\\\vdots\\\Lat_{N-1,N}}
(\Lat_{N,N})^{-1}
\mat{c}{(\Lat_{N,1})^*\\\vdots\\(\Lat_{N,N-1})^*}^*.
\end{equation*}
Furthermore, if $\sU$ is f{}inite dimensional, then all computations
involve f{}inite matrices only. In this sense the relaxed Nehari problem
is very different from the classical Nehari problem.

\noindent{\bf Proof of Theorem \ref{th:Neh}.}
Since $R$ in
(\ref{RandQ}) is an isometry, we have that the lifting data set
$\Om=\LDS$, with $U'$ the Sz.-Nagy-Sch\"affer isometric lifting of
$T'$ in (\ref{defTpr}) and $R$ and $Q$ as in (\ref{RandQ}), has the
property that $A$ is a strict contraction (by assumption) and $R$ is
left invertible. So the description of all contractive interpolants
in Theorem \ref{mainth} can be applied to this particular lifting
data set. Since also $Q$ is an isometry, we have $R^*R=Q^*Q$. Hence
$D_\circ=0$ on $\sD_\circ=\{0\}$, where $D_\circ$ and $\sD_\circ$
are given by (\ref{Dcirc}). We already made the identif{}ication
between $\sD_{T'}$ and $\sY$. Moreover, we have $\kr
Q^*=\sU\oplus\{0\}^{N-1}\subset\sU^N$ and $\kr
R^*=\{0\}^{N-1}\oplus\sU\subset\sU^N$. So both $\kr Q^*$ and $\kr
R^*$ can be identif{}ied with $\sU$. In that case the projections
$\Pi_{\kr Q^*}$ and $\Pi_{\kr R^*}$ are given by
\begin{equation}\label{kerns}
\Pi_{\kr Q^*}=\mat{cccc}{I_\sU&0&\cdots&0}:\sU^N\to\sU\ands \Pi_{\kr
R^*}=\mat{cccc}{0&\cdots&0&I_\sU}:\sU^N\to\sU.
\end{equation}
In particular, $\Pi_{\kr Q^*}=E_N^*$, with $E_N$ the operator in
(\ref{NehCons}).

Now let $H$ be a function from $\eH^2(\sU,\sY)$ with Taylor
coeff{}icients $(H_n)_{n\in\BN}$ at zero. Then the operator $\Ga_H$
from $\sU$ into $H^2(\sY)$ def{}ined in (\ref{multH2}) is given by
\[
\Ga_H=\Four_\sY\Pi_{\ell^2_+(\sY)}\tilB E_N=\Pi_{H^2(\sY)}\Psi\tilB
E_N,
\]
where $\tilB$ is given by (\ref{defB}), and $\Four_\sY$ and $\Psi$
are the unitary operators def{}ined in (\ref{four}) and
(\ref{uneqlifts}), respectively. According to Proposition
\ref{pr:Neh} and the last but one paragraph of the previous
subsection we have that $(H_n)_{n\in\BN}$ is an $N$-complementary
sequence associated with $(F_{-n-1})_{n\in\BN}$ if and only if
$\Psi\tilB$ is a contractive interpolant for $\Om$. Moreover, all
contractive interpolants for $\Om$ are obtained in this way. Then
Theorem \ref{mainth} and the above analysis of the spaces $\kr Q^*$,
$\kr R^*$, $\sD_\circ$ and $\sD_{T'}$ imply that $(H_n)_{n\in\BN}$
is an $N$-complementary sequence if and only if there exists a Schur
class function $\Va$ from $\eS(\sU,\sY\oplus\sU)$ such that
\begin{equation*}
H(\la)=\Pa_{22}(\la)E_N+\Pa_{21}(\la)\Va(\la)
(I-\Pa_{11}(\la)\Va(\la))^{-1}\Pa_{12}(\la)E_N\quad(\la\in\BD),
\end{equation*}
where $\Pa_{11}$, $\Pa_{12}$, $\Pa_{21}$ and $\Pa_{22}$ are the
functions def{}ined in Theorem \ref{mainth}. Put for each $\la\in\BD$
\begin{equation}\label{Pc`s}
\Pc_{11}(\la)=\Pa_{11}(\la),\quad
\Pc_{12}(\la)=\Pa_{12}(\la)E_N,\quad
\Pc_{21}(\la)=\Pa_{21}(\la)\ands
\Pc_{22}(\la)=\Pa_{22}(\la)E_N.
\end{equation}

Then we obtain from Theorem \ref{mainth} that
$\Pc_{11}\in\eS(\sY\oplus\sU,\sU)$, $\Pc_{12}\in\eH^2(\sU,\sU)$,
$\Pc_{21}\in\eS(\sY\oplus\sU,\sY)$ and $\Pc_{22}\in\eH^2(\sU,\sY)$.
Note that $\Ga_{\Pc_{12}}=\Ga_{\Pa_{12}}E_N$,
$\Ga_{\Pc_{22}}=\Ga_{\Pa_{22}}E_N$ and $\tilGa_-=AE_N$.
Hence the operator $\hat M$ in (\ref{hatM}) is obtained from $M$ in
(\ref{multop}) after multiplication from the right by
\[
\mat{cc}{I_{H^2(\sY\oplus\sU)}&0\\0&E_N}
:\mat{c}{H^2(\sY\oplus\sU)\\\sU}\to\mat{c}{H^2(\sY\oplus\sU)\\\sU^N}.
\]
Since $R^*R=Q^*Q$, we obtain from Corollary \ref{cor:multop} that $M$,
and thus $\hat M$, is an isometry if $X_1$ in (\ref{X1X5}) is pointwise
stable. This proves to be the case. In fact, we have that $X_1$ in
(\ref{X1X5}) has $\spec(X_1)<1$. To see this it suff{}ices to show
that $R-\la Q$ is left invertible for each $\la$ in $\BD$. See
Proposition 5.3 in \cite{FFK02a}, and the remark in the third
paragraph after Theorem \ref{mainth}. Indeed, we have that $R-\la Q$
is left invertible for each $\la\in\BD$. In fact, a left inverse is
given by
\[
\mat{ccccc}
{I&0&\cdots&\cdots&0\\
\la I&I&\ddots&&\vdots\\
\vdots&\ddots&\ddots&0&\vdots\\
\la^{N-2}I&\cdots&\la I&I&0}:\sU^N\to\sU^{N-1}.
\]
So to complete the proof it remains to show that the functions
$\Pc_{11}$, $\Pc_{12}$, $\Pc_{21}$ and $\Pc_{22}$ def{}ined in
(\ref{Pc`s}) can also be written as in (\ref{NehPhi2}), and that
$\spec(\TA)<1$. The later immediately shows that $\Pc_{12}$ and
$\Pc_{22}$ are functions in $\eH^\infty(\sU,\sU)$ and
$\eH^\infty(\sU,\sY)$, respectively.

{}From (\ref{kerns}), the def{}initions of $\De_Q$ and $\De_R$ in
(\ref{Deltas}) and the assumption that $D_A^{-2}$ is given by
(\ref{defectAinv}) we immediately obtain that $\De_Q=\Lat_{1,1}$
and $\De_R=\Lat_{N,N}$. Moreover, the $N-1\times N-1$ left upper
corner and $N-1\times N-1$ right lower corner of $D_A^2$ in
(\ref{defectA}) are equal to $R^*D_A^2R$ and $Q^*D_AQ$, respectively.
The fact that $D_{T'}$ is given by (\ref{defectT'}) and the def{}inition
of $A$ in (\ref{NHankel}) show that $D_{T'}A=F$, where $F$ is given by
(\ref{NehCons}). Thus $D_{T'}AR=FR=\mat{ccc}{Y_{-1}&\cdots&Y_{-N+1}}$.
This implies that $G$ in (\ref{NehCons}) can also be written as
$G^*=R(R^*D_A^2R)^{-1}R^*A^*D_{T'}$. Using the fact that $D_\circ=0$ on
$\sD_\circ=\{0\}$, we obtain that
\[
\De_\Om=I+D_{T'}AR(R^*D_A^2R)^{-1}R^*A^*D_{T'}=I+FG^*.
\]
So we have the following identities:
\begin{equation*}\label{Nidens}
\De_Q=\Lat_{1,1},\quad
\De_R=\Lat_{N,N},\quad
\De_\Om=I+FG^*,\quad
F=D_{T'}A\ands
G=D_{T'}AR(R^*D_A^2R)^{-1}R^*.
\end{equation*}

Let $X_1$, $X_2$, $X_3$, $X_4$ and $X_5$ be the operators def{}ined in
(\ref{X1X5}). Then we immediately obtain that
\begin{equation}\label{NX3X4X5}
\begin{array}{l}
X_3=(\Lat_{1,1})^{-\half}E_N^*:\sU^N\to\sU,\quad
X_4=FX_1:\sU^N\to\sY\\[.2cm]
X_5=(I+FG^*)^{-\half}\Pi_{\sY}:\sY\oplus\sU\to\sY.
\end{array}
\end{equation}
Writing out $\Pc_{11}$, $\Pc_{12}$, $\Pc_{21}$ and $\Pc_{22}$ in terms
of $X_1$, $X_2$, $X_3$, $X_4$ and $X_5$ with the identities for $X_3$,
$X_4$ and $X_5$ in (\ref{NX3X4X5}) we obtain that
\begin{equation}\label{Pcinterforms}
\begin{array}{rcl}
\Pc_{11}(\la)&=&\la X_3(I-\la X_1)^{-1}X_2
=\la(\Lat_{1,1})^{-\half}E_N^*(I-\la X_1)^{-1}X_2,\\[.2cm]
\Pc_{12}(\la)&=&X_3E_N+\la X_3(I-\la X_1)^{-1}X_1E_N\\[.2cm]
&=&(\Lat_{1,1})^{-\half}
+\la(\Lat_{1,1})^{-\half}E_N^*(I-\la X_1)^{-1}X_1E_N,\\[.2cm]
\Pc_{21}(\la)&=&X_5+\la X_4(I-\la X_1)^{-1}X_2
=X_5+\la FX_1(I-\la X_1)^{-1}X_2\\[.2cm]
&=&X_5-FX_2+F(I-\la X_1)^{-1}X_2,\\[.2cm]
\Pc_{22}(\la)&=&X_4(I-\la X_1)^{-1}E_N=F(I-\la X_1)^{-1}
X_1E_N.
\end{array}
\quad\quad(\la\in\BD)
\end{equation}
Therefore, to complete the proof, it suff{}ices to show that
\begin{equation}\label{NX1X2}
X_1=\TA\ands
X_2=-\mat{cc}{G^*(I+FG^*)^{-\half}&C_2}.
\end{equation}
Indeed, assume that (\ref{NX1X2}) holds. Then $C_1=-X_1E_N$ and
\begin{eqnarray*}
X_5-FX_2&=&\mat{cc}{(I+FG^*)^{-\half}&0}
+\mat{cc}{FG^*(I+FG^*)^{-\half}&FC_2}\\
&=&\mat{cc}{(I+FG^*)(I+FG^*)^{-\half}&FC_2}
=\mat{cc}{(I+FG^*)^\half&FC_2}.
\end{eqnarray*}
This computation combined with the formulas for $\Pc_{11}$,
$\Pc_{12}$, $\Pc_{21}$ and $\Pc_{22}$ found in (\ref{Pcinterforms})
show that (\ref{NehPhi2}) holds, and also that
$\spec(\TA)=\spec(X_1)<1$.

We will now prove (\ref{NX1X2}), starting with the f{}irst identity.
First we claim that $X_1$ can be written as
\[
X_1=RQ^*+R(Q^*D_A^2Q)^{-1}Q^*D_A^2P_{\kr Q^*}.
\]
To see that this is the case, observe that both $Q^*$ and
$(Q^*D_A^2Q)^{-1}Q^*D_A^2$  are left inverses of $Q$. Certainly all
left inverses of $Q$ are equal on $\im Q$. Since
$Q^*|(\sU^N\ominus\im Q)=Q^*|\kr Q^*=0$, we obtain that
$(Q^*D_A^2Q)^{-1}Q^*D_A^2=Q^*+(Q^*D_A^2Q)^{-1}Q^*D_A^2P_{\kr Q^*}$,
which proves our claim.

Note that $RQ^*$ is the backward shift on
$\sU^N$. This implies that $\TA=RQ^*+C_1E_N^*$. Using that
$P_{\kr Q^*}=E_NE_N^*$ and $\Pi_{\kr R^*}C_1=0$ we obtain that
the f{}irst identity in (\ref{NX1X2}) holds if
$R^*C_1=-(Q^*D_AQ)^{-1}Q^*D_A^2E_N$. This identity follows from
the next computation:
\begin{eqnarray*}
Q^*D_A^2QR^*C_1\Lat_{1,1}
&=&Q^*D_A^2Q\mat{c}{\Lat_{2,1}\\\vdots\\\Lat_{N,1}}
=Q^*D_A^2\mat{c}{0\\\Lat_{2,1}\\\vdots\\\Lat_{N,1}}\\
&=&Q^*D_A^2D_A^{-2}E_N-Q^*D_A^2E_N\Lat_{1,1}\\
&=&Q^*E_N-Q^*D_A^2E_N\Lat_{1,1}
=-Q^*D_A^2E_N\Lat_{1,1}.
\end{eqnarray*}
Indeed, to obtain the desired equality multiply the f{}irst and the last
term in the above sequence of equalities with $(Q^*D_A^2Q)^{-1}$ from
the left and with $(\Lat_{1,1})^{-1}$ from the right.

To see that the identity for $X_2$ in (\ref{NX1X2}) holds, note that
$X_2$ can be written as
\[
X_2=-\mat{cc}{R(R^*D_A^2R)^{-1}R^*A^*D_{T'}\De_\Om^{-\half}
&D_A^{-2}\Pi_{\kr R^*}\De_R^{-\half}}:\mat{c}{\sY\\\sU}\to\sU^N.
\]
Earlier we obtained that $\De_\Om=I+FG^*$ and
$R(R^*D_A^2R)^{-1}R^*A^*D_{T'}=G^*$. This shows that
$X_2|\sY=-G^*(I+FG^*)^{-\half}$. The fact that $X_2|\sU=-C_2$ follows
from the def{}inition of $\Pi_{\kr R^*}$ in (\ref{kerns}),
the formula for $D_A^{-2}$ in (\ref{defectAinv}) and the identity
$\De_R=\Lat_{N,N}$. Thus the second identity in (\ref{NX1X2})
holds.
\epr

%%%%%%%%%%%%%%%%%%%%%%%%%%%%%%%%%%%%%%%%%%%%%%%%%%%%%%%%%%%%%%%%%%%%%%%%
\subsection{Special cases}

As an illustration we specify the result in Theorem \ref{th:Neh} for
two special cases, namely, the case that $N=1$ and the case that
$F_n=0$ for $n=-1,-2,\ldots.$

%%%%%%%%%%%%%%%%%%%%%%%%%%%%%%%%%%%
\begin{corollary}\label{cor:N=1}
Let $F_{-1},F_{-2},\ldots$ be operators in $\Ops(\sU,\sY)$ such that
$\sum_{n=1}^\infty\|F_{-n}u\|^2<\infty$ for each $u\in\sU$.
Assume that the 1-truncated Hankel operator $A$ is a contraction. Let
$H$ be a function from $\eH^2(\sU,\sY)$. Then the Taylor coeff{}icients
$(H_n)_{n\in\BN}$ of $H$ at zero form a $1$-complementary sequence
associated with  $(F_{-n-1})_{n\in\BN}$ if and only if there exists a
Schur class function $V$ from $\eS(\sD_A,\sY\oplus\sD_A)$ such that
\begin{equation}\label{N=1}
H(\la)=\Pi_\sY V(\la)(I-\la\Pi_{\sD_A}V(\la))^{-1}D_A
\quad\quad(\la\in\BD).
\end{equation}
Moreover, if $A$ is a strict contraction, then $\sD_A=\sU$.
\end{corollary}

Note that in Corollary \ref{cor:N=1} it is not required
that the 1-truncated Hankel operator is a strict contraction. We
give two different proofs. The f{}irst is based on Theorem \ref{th:FtHK}
and Proposition \ref{pr:Neh}. In the second proof, assuming $\|A\|<1$,
we show that Theorem \ref{th:Neh} reduces to Corollary \ref{cor:N=1},
when specif{}ied for $N=1$.

\noindent{\bf Proof 1} (Using Theorem \ref{th:FtHK} and
Proposition \ref{pr:Neh}).
Since $N=1$, we see that
$A=\mat{ccc}{\cdots&F_{-2}^*&F_{-1}^*}^*$ from $\sU$ to
$\ell^2_-(\sY)$, and that $R$ and $Q$ in (\ref{RandQ}) reduce to
$R=Q=0$ from $\{0\}$ to $\sU$. Note that $E_N=I_\sU$, with $E_N$ the
operator def{}ined in Proposition \ref{pr:Neh}. Then Proposition
\ref{pr:Neh} shows that the Fourier coeff{}icients of $H$ form a
$1$-complementary sequence associated with $(F_{-n-1})_{n\in\BN}$ if
and only if $\Ga_H=\Pi_{H^2(\sY)}B$ for a contractive interpolant
$B$ for $\LDS$, where $U'$ is the Sz.-Nagy-Sch\"affer isometric
lifting of the contraction $T'$ def{}ined in (\ref{defTpr}), see also
the second paragraph of the proof of Theorem \ref{th:Neh}.

Since $R=Q=0$, we obtain that $\om=0$ and $\sF=\{0\}$, see
(\ref{om}). Therefore any function $Z$ from
$\eS(\sD_A,\sY\oplus\sD_A)$ has the property that $Z(\la)|\sF=\om$
for each $\la\in\BD$. Recall that $\sD_{T'}=\sY$, see
(\ref{defectT'}). Applying Theorem \ref{th:FtHK}, we see that the
Fourier coeff{}icients of $H$ form a $1$-complementary sequence
associated with $(F_{-n-1})_{n\in\BN}$ if and only if $H$ is given
by (\ref{N=1}), where $V$ is from the Schur class
$\eS(\sD_A,\sY\oplus\sD_A)$. If $A$ is a strict contraction, then
$D_A$ is invertible on $\sU$. In particular, this implies that
$\sD_A=\sU$.
\epr

\noindent{\bf Proof 2} (Using Theorem \ref{th:Neh} and assuming that
$\|A\|<1$).
Assume that the $1$-truncated Hankel operator for
$(F_{-n-1})_{n\in\BN}$ is a strict contraction.
Let $\Pc_{11}$, $\Pc_{12}$, $\Pc_{21}$ and $\Pc_{22}$ be the functions
def{}ined in Theorem \ref{th:Neh} for the case $N=1$. Then the operators
$C_1$, $C_2$, $F$ and $G$ in (\ref{NehCons}) reduce to
\[
C_1=0\ons\sU,\quad
C_2=(\Lat_{1,1})^\half\ons\sU,\quad
F=F_{-1}:\sY\to\sU\ands
G=0:\sY\to\sU.
\]
This implies that
\[
\TA=0\ons\sU,\quad E_1=I\ons\sU\ands
\mat{cc}{G^*(I+FG^*)&C_2}=\mat{cc}{0&(\Lat_{1,1})^\half}.
\]
Moreover, observe that $(\Lat_{1,1})^\half=D_A^{-1}$.
Therefore we have for each $\la\in\BD$ that
\begin{eqnarray*}
\Pc_{11}(\la)&=&-\la D_A(I-\la0)^{-1}\mat{cc}{0&D_A^{-1}}
=-\mat{cc}{0&I}=-\Pi_\sY,\\
\Pc_{12}(\la)&=&D_A-\la D_A(I-\la0)^{-1}0=D_A,\\
\Pc_{21}(\la)&=&\mat{cc}{(I+0)^{\half}&F_{-1}D_A^{-1}}
-F_{-1}(I-\la0)^{-1}\mat{cc}{0&D_A^{-1}}=\mat{cc}{I&0}=\Pi_\sU,\\
\Pc_{22}(\la)&=&-F_1(I-\la0)^{-1}0=0.
\end{eqnarray*}

According to Theorem \ref{th:Neh} the Fourier coeff{}icients of $H$ form
an $N$-complementary sequence associated with $(F_{-n-1})_{n\in\BN}$
if and only if
\[
H(\la)=\Pi_\sY V(\la)(I+\la\Pi_\sU V(\la))^{-1}D_A\quad(\la\in\BD),
\]
for some $V\in\eS(\sU,\sY\oplus\sU)$. Note that instead of the plus
in $(I+\Pi_\sU V(\la))$ we can also take a minus since
\[
V\in\eS(\sU,\sY\oplus\sU)\quad\mbox{if and only if}\quad
\mat{c}{\Pi_\sY V(\cdot)\\-\Pi_{\sU}V(\cdot)}
\in\eS(\sU,\sY\oplus\sU).
\]
\epr

Finally, we consider the relaxed Nehari problem with
$F_{-1}=F_{-2}=\cdots=0$. That is, we seek a description of all
functions $H$ from $\eH^2(\sU,\sY)$ with Fourier coeff{}icients
$(H_n)_{n\in\BN}$ at zero such that the operator
\begin{equation}\label{L0}
\mat{ccccc}{H_0&0&0&\cdots&0\\H_1&H_0&0&\cdots&0\\
H_2&H_1&H_0&\cdots&0\\\vdots&\vdots&\vdots&\ddots&\vdots\\
H_{N-1}&H_{N-2}&H_{N-3}&\cdots&H_0\\
\vdots&\vdots&\vdots&\vdots&\vdots}:\sU^N\to\ell^2_+(\sY)
\end{equation}
is a contraction. By associating with each function $H$ from
$\eH^2(\sU,\sY)$ the norm of the operator in (\ref{L0}) we induce a
Banach space structure on the set $\eH^2(\sU,\sY)$. This Banach space
appears in \cite{FFK02a}, in the contexts of certain interpolation
problems, and for the case that $\sU$ and $\sY$ are f{}inite
dimensional, in \cite{ABL96} and \cite{ABL97}.

By specifying Theorem \ref{th:Neh} for the case that
$F_{-1}=F_{-2}=\cdots=0$ we obtain the following description of all
functions $H$ from $\eH^2(\sU,\sY)$ with $L_0$ in (\ref{L0})
contractive.

%%%%%%%%%%%%%%%%%%%%%%%%%%%%%%%%%%%
\begin{corollary}\label{cor:Fn=0}
Let $H$ be a function from $\eH^2(\sU,\sY)$ with Fourier coeff{}icients
$(H_n)_{n\in\BN}$ at zero. Then the operator in \tu{(\ref{L0})} is a
contraction if and only if there exists a Schur class function $V$
from $\eS(\sU,\sY\oplus\sU)$ such that
\begin{equation}\label{Fn=0}
H(\la)=\Pi_\sY V(\la)(I-\la^N\Pi_\sU V(\la))^{-1}\quad(\la\in\BD).
\end{equation}
\end{corollary}

Corollary \ref{cor:Fn=0} is an operator-valued version of a result from
\cite{ABL96}. The case $N=1$ appears as a corollary in \cite{FtHK06a}
and is fundamental in the proof of the f{}irst main result in
\cite{FtHK06b}.

\noindent{\bf Proof of Corollary \ref{cor:Fn=0}.}
Note that the $N$-truncated Hankel operator $A$ for the sequence
$F_{-1}=F_{-2}=\cdots=0$ is the zero operator from $\sU^N$ to
$\ell^2_-(\sY)$. In particular, $A$ is a strict contraction.
So we can apply the result of Theorem \ref{th:Neh} to this
Nehari data.

Note that in this case both $D_A^2$ in (\ref{defectA})
and $D_A^{-2}$ in (\ref{defectAinv}) are equal to the identity
operator on $\sU^N$. In particular, $\Lat_{i,j}=I$ if $i=j$ and
$\Lat_{i,j}=0$ if $i\not=j$, for $i,j=1,\ldots,N$. Moreover, the
operators $G_1,\ldots,G_{N-1}$ in (\ref{G1GN-1}) are all equal to
$0$. Therefore we obtain that the operators $F$, $G$, $C_1$ and $C_2$ in (\ref{NehCons}) are given by $F=G=0$ from $\sU^N$ to $\sY$, $C_1=0$ from
$\sU$ to $\sU^N$ and $C_2^*=\mat{cccc}{0&\cdots&0&I}$ from $\sU^N$
to $\sU$. Furthermore, $\TA$ in (\ref{NehCons}) reduces to
\begin{equation*}
\TA=\mat{ccccc}{0&I&0&\cdots&0\\
\vdots&\ddots&\ddots&\ddots&\vdots\\
\vdots&&\ddots&\ddots&0\\
\vdots&&&\ddots&I\\
0&\cdots&\cdots&\cdots&0},
\quad\mbox{so}\quad
(I-\la\TA)^{-1}
=\mat{cccc}{I&\la I&\cdots&\la^{N-1}I\\0&I&\ddots&\vdots\\
\vdots&\ddots&\ddots&\la I\\0&\cdots&0&I}
\end{equation*}
for each $\la\in\BD$. Using the form of $C_2$ found above, and the
def{}inition of $E_N$ in (\ref{NehCons}), we see that
$E_N^*(I-\la\TA)^{-1}C_2=\la^{N-1}I$ for each $\la\in\BD$.
Therefore we have for each $\la\in\BD$ that
\begin{eqnarray*}
\Pc_{11}(\la)
&=&-\la I^{-\half}E_N^*(I-\la\TA)^{-1}\mat{cc}{0&C_2}
=-\mat{cc}{0&\la^NI}=-\la^N\Pi_\sU\\
\Pc_{12}(\la)
&=&I^{-\half}-\la I^{-\half}E_N^*(I-\la\TA)^{-1}0=I\\
\Pc_{21}(\la)
&=&\mat{cc}{(I+0)^\half&0C_2}-0(I-\la\TA)^{-1}\mat{cc}{0&C_2}
=\mat{cc}{I&0}=\Pi_\sY\\
\Pc_{22}(\la)
&=&-0(I-\la\TA)^{-1}0=0.
\end{eqnarray*}
Inserting these formulas for $\Pc_{11}$, $\Pc_{12}$, $\Pc_{21}$ and
$\Pc_{22}$ into (\ref{NehPhi1}) we obtain (\ref{Fn=0}).
\epr

\noindent{\bf Acknowledgement.}
 The author thanks Prof. M.A.
Kaashoek and Prof. A.E. Frazho for their many useful comments and
suggestions.

%%%%%%%%%%%%%%%%%%%%%%%%%%%%%%%%%%%%%%%%%%%%%%%%%%%%%%%%%%%%%%%%%%%%%%%%
%%%%%%%%%%%%%%%%%%%%%%%%%%%%%%%%%%%%%%%%%%%%%%%%%%%%%%%%%%%%%%%%%%%%%%%%

\end{document}